\documentclass[11pt, reqno]{amsart}
\usepackage{indentfirst, amssymb, amsmath, amsthm, mathrsfs, setspace, indentfirst, enumerate,  mathrsfs, amsmath, amsthm}
\usepackage[bookmarksnumbered, colorlinks, plainpages]{hyperref}
\usepackage{mathrsfs}

\newtheorem*{theoA}{Theorem A}
\newtheorem*{theoB}{Theorem B}
\newtheorem*{theoC}{Theorem C}

\newtheorem{theo}{Theorem}[section]
\newtheorem{lem}{Lemma}[section]
\newtheorem{cor}{Corollary}[section]

\newtheorem{exm}{Example}[section]

\newcommand{\ol}{\overline}
\newcommand{\ul}{\underline}

\newcommand{\be}{\begin{equation}}
\newcommand{\ee}{\end{equation}}
\newcommand{\beas}{\begin{eqnarray*}}
\newcommand{\eeas}{\end{eqnarray*}}
\newcommand{\bea}{\begin{eqnarray}}
\newcommand{\eea}{\end{eqnarray}}

\numberwithin{equation}{section}
\begin{document}

\title[E\MakeLowercase{ntire solutions of} F\MakeLowercase{ermat-Type Partial} ......]{
\LARGE E\Large\MakeLowercase {ntire solutions of} F\MakeLowercase {ermat}-\MakeLowercase {Type Partial Differential-Difference equations in higher dimensions}}

\date{}
\author[S. M\MakeLowercase {ajumder}, J. B\MakeLowercase{anerjee}, A. B\MakeLowercase{anerjee} \MakeLowercase{and} E. G. S\MakeLowercase{aleeby}]
{S\MakeLowercase {ujoy} M\MakeLowercase {ajumder}$^1$, J\MakeLowercase {hilik} B\MakeLowercase {anerjee}$^2$, A\MakeLowercase{bhijit} B\MakeLowercase {anerjee}$^3$ 
\MakeLowercase{and} E\MakeLowercase {lias} S\MakeLowercase {aleeby}$^4$$^*$}

\address{$^1$ Sujoy Majumder: Department of Mathematics, Raiganj University, Raiganj, West Bengal 733134, India}
\email{sm05math@gmail.com, sjm@raiganjuniversity.ac.in}

\address{$^2$ Jhilik Banerjee: Department of Mathematics, University of Kalyani, West Bengal 741235, India}
\email{jhilikbanerjee38@gmail.com, jhilikmath24@klyuniv.ac.in}

\address{$^3$ Abhijit Banerjee: Department of Mathematics, University of Kalyani, West Bengal 741235, India}
\email{abanerjee\_kal@yahoo.co.in}

\address{$^4$ Elias G. Saleeby: Department of Mathematics, Al Akhawayn University, Ifrane 53000, Morocco}
\email{esaleeby@yahoo.com, e.saleeby@aui.ma}

\thanks{*\emph{Corresponding Author}: Elias G. Saleeby.}

\begin{abstract} The purpose of this paper is to investigate the finite-order entire solutions of the following Fermat-type partial differential-difference equations:
\begin{align*} (P_1(z)f(z+c))^{m_1}+\left(P_2(z)f(z)+\sideset{}{_{j=1}^m}{\sum} Q_j(z)\partial_{z_j}(f(z))\right)^{m_2}=P(z),
\end{align*}
\begin{align*} (P_1(z)f(z))^{m_1}+\left(P_2(z)f(z+c)+\sideset{}{_{j=1}^m}{\sum} Q_j(z)\partial_{z_j}(f(z))\right)^{m_2}=P(z)
\end{align*}
and 
\begin{align*} \left(\sideset{}{_{j=1}^m}{\sum} Q_j(z)\partial_{z_j}(f(z))\right)^{m_1}+\left(P_1(z)f(z+c)+P_2(z)f(z)\right)^{m_2}=P(z),
\end{align*}
in $\mathbb{C}^m$, where $m_1,m_2\in\mathbb{N}$ such that $m_1+m_2>2$, $c=(c_1,c_2,\ldots,c_m)\in\mathbb{C}^m$, $P(\not\equiv 0)$, $P_1$, $P_2$ and $Q_i\;(i=1,2,\ldots,m)$ are polynomials in $\mathbb{C}^m$ such that $P(z)$ is irreducible in $\mathbb{C}^m$. 
 The results established in this paper generalize those of Xu et al. \cite{Xu-Luo-Cao-2025} from $\mathbb{C}$ to $\mathbb{C}^m$.
\end{abstract}

\thanks{Typeset by \AmS -\LaTeX}
\maketitle

\section{{\bf Introduction}}

\noindent
Fermat's Last Theorem states that the equation $x^n + y^n = 1$ has no rational solutions with $xy\neq 0$ for $n\geq 3$ (see \cite{TW1, W1}), whereas the equation $x^2 + y^2 = 1$ does admit nontrivial rational solutions. Analogously, the functional equation
\begin{align}\label{ne1}
f^n(z) + g^n(z) = 1,
\end{align}
where $n$ is a positive integer, can be regarded as the Fermat Diophantine equation $x^n + y^n = 1$ over function fields.
In 1927, Montel \cite{M1} proved that Eq. \eqref{ne1} has no transcendental entire solutions for $n \geq 3$. Furthermore, as shown in \cite[Theorem 4]{FG2}, Eq. \eqref{ne1} actually has no non-constant entire solutions when $n > 2$. Gross \cite[Theorem 3]{FG1} also proved that Eq. \eqref{ne1} has no non-constant meromorphic solutions when $n > 3$.
For $n = 2$, Gross \cite[Theorem 4]{FG2} found that Eq. \eqref{ne1} possesses entire solutions of the form $f(z) = \sin(h(z))$ and $g(z) = \cos(h(z))$, where $h(z)$ is an entire function. For $n = 2$ and $n = 3$, Gross \cite[Theorem 1]{FG1} and Baker \cite[Theorem 1]{Baker}, respectively, obtained the solutions of Eq. \eqref{ne1} and its specific forms.

\noindent
In 1970, Yang \cite{Y1} investigated the generalized Fermat-type equation
\begin{align*}
f^{m_1}(z) + g^{m_2}(z) = 1
\end{align*}
and demonstrated that it has no non-constant entire solutions when $\frac{1}{m_1} + \frac{1}{m_2} < 1$ (see the proof of Theorem 1). Consequently, the equation has no non-constant entire solutions when $m_1 > 2$ and $m_2 > 2$. However, for the case where $m_1 = m_2 = 2$ and $g(z)$ has a specific relationship with $f(z)$, many researchers have investigated the existence of solutions, leading to a substantial body of published work (see \cite{Liu1}-\cite{Lu2}, \cite{TL1}, \cite{YL1}), and research in this area remains active.
\medskip

\noindent
Throughout the paper, we use the notations $\ol{\omega}(z)=\omega(z+c)$ and $\ul{\omega}(z)=\omega(z-c)$ for any meromorphic function $\omega(z)$ and a finite nonzero complex number $c$.
Recently, Xu et al. \cite{Xu-Luo-Cao-2025} investigated the entire solutions $f(z)$ of Fermat-type difference-differential equations involving both the derivative operator $f'(z)$ and the difference operator $f(z + c)$. Specifically, they determined all entire solutions $f$ of the functional equation
$$X^2 + Y^2 = p(z),$$
satisfying one of the following three cases:
\begin{enumerate}
	\item[(i)] $X = p_1(z)f(z + c)$ and $Y = p_2(z)f(z) + p_3(z)f'(z)$;
	\item[(ii)] $X = p_1(z)f(z)$ and $Y = p_2(z)f'(z) + p_3(z)f(z + c)$;
	\item[(iii)] $X = p_1(z)f'(z)$ and $Y = p_2(z)f(z + c) + p_3(z)f(z)$,
\end{enumerate}
where $c$ is a nonzero complex number, $p_1$, $p_2$, and $p_3$ are polynomials in $\mathbb{C}$ satisfying $p_1 p_3 \not\equiv 0$, and $p$ is a nonzero irreducible polynomial in $\mathbb{C}$. 
Xu et al. \cite{Xu-Luo-Cao-2025} proved the following results:
\begin{theoA} \cite[Theorem 1.1]{Xu-Luo-Cao-2025}
Let $f$ be an entire solution of the Fermat-type difference-differential equation
\begin{align*}
\left[p_1(z)f(z+c)\right]^2
+\left[p_2(z)f(z)+p_3(z)f'(z)\right]^2
=p(z),
\end{align*}
where the polynomials $p_1,p_2$ and $p_3$ in $\mathbb{C}$ satisfy
$p_1p_3\not\equiv 0$, and $p$ is a nonzero irreducible polynomial
in $\mathbb{C}$. Then $f$ has the form
\begin{align*}f(z)=\frac{\exp{(\iota g(z))}+p \exp{(-\iota g(z))}}{2p_1(z)},\end{align*}
where $g$ is a constant or non-constant linear polynomial
$g(z)=Az+B$ such that
\begin{align*}
\exp{(-\iota Ac)}\equiv\frac{k p_1^2}{(p_2-Ap_3)p_1-\iota p_3p_1'},\qquad (k=\pm 1).
\end{align*}
\end{theoA}

\begin{theoB} \cite[Theorem 1.3]{Xu-Luo-Cao-2025}	
Let $f$ be an entire solution of the Fermat-type difference-differential equation
\begin{align*}
\left[p_1(z)f(z)\right]^2+\left[p_2(z)f'(z)+p_3(z)f(z+c)\right]^2=p(z),
\end{align*}
where the polynomials $p_1$, $p_2$ and $p_3$ in $\mathbb{C}$ satisfy
$p_1p_3\not\equiv 0$, and $p$ is a nonzero irreducible polynomial
in $\mathbb{C}$. Then $f$ has the form
\begin{align*}f(z)=\frac{\exp{(\iota g(z))}+p(z)\exp{(-\iota g(z))}}{2p_1(z)},\end{align*}
where $g$ is a constant or non-constant linear polynomial
$g(z)=Az+B$ such that
\begin{align*}\exp{(-\iota Ac)}\equiv\frac{\iota p_1^2p_3}{p_1(Ap_2p_1+\iota p_2'+kp_1^2)},\qquad (k=\pm 1).\end{align*}
\end{theoB}

\begin{theoC} \cite[Theorem 1.5]{Xu-Luo-Cao-2025}	
Let $p_1$, $p_2$ and $p_3$ be three polynomials in $\mathbb{C}$
satisfying $p_1p_3\not\equiv 0$ and $p_2=Mp_3$,
where $M$ is a constant, and let $p$ be a nonzero irreducible
polynomial in $\mathbb{C}$. If $f$ is an entire solution of the
Fermat-type difference-differential equation
\begin{align*}
\left[p_1(z)f'(z)\right]^2+\left[p_2(z)f(z+c)+p_3(z)f(z)\right]^2
=p(z),
\end{align*}
then $f$ satisfies
\begin{align*}f'(z)=\frac{\exp{(\iota g(z))}+p(z)\exp{(-\iota g(z))}}{2p_1(z)},\end{align*}
where $g$ satisfies the following two cases:
\begin{enumerate}
\item[\em{(i)}] Whenever $M=0$ (that is, $p_2\equiv 0$), then $g$
is a constant or non-constant linear polynomial
\begin{align*}g=\frac{1}{k}\frac{p_3}{p_1}z+\text{Constant},\end{align*}
such that $p$, $p_1$ and $p_3$ are constant.	
\item[\em{(ii)}] Whenever $M\neq 0$, then $g$ is a constant or
non-constant linear polynomial $g(z)=Az+B$ such that
\begin{align*}\exp{(-\iota Ac)}\equiv\frac{\iota p_1p_2p_3}{k \iota A p_1p_3\overline{p_1}
-\iota^2p_3^2\overline{p_1}-k p_3^3p_1\overline{p_1}}.\end{align*}
\end{enumerate}
\end{theoC}
\bigskip

\noindent
In recent years, the classical framework in $\mathbb{C}$ has been broadly extended to higher dimensions. Fermat-type partial differential-difference equations in several complex variables ($\mathbb{C}^m$) combine differential operators (such as partial derivatives $\frac{\partial f}{\partial z_j}$) with shift/difference operators ($f(z+c)$). These equations are generally expressed in the form:
\begin{align*}
\left(L_1(f)\right)^{m_1} + \left(L_2(f)\right)^{m_2} = P(z),
\end{align*}
where $L_1$ and $L_2$ are linear differential, difference, or hybrid partial differential-difference operators with polynomial coefficients, $P(z)$ is a given polynomial or entire function, and $m_1, m_2 \in \mathbb{N}$.

\noindent
Investigating solutions to these equations in $\mathbb{C}^m$ relies heavily on several complex variables Nevanlinna theory, including growth characteristics, characteristics of logarithmic derivatives, and value distribution mechanisms across complex hyperplanes. The primary objectives in this area are:
\begin{enumerate}
\item[1.] Solvability Conditions: Establishing strict constraints on the exponents $m_1$ and $m_2$. The condition $m_1 + m_2 > 2$ only excludes the trivial case.
\item[2.] Explicit Factorization \& Structural Representation: Characterizing the precise exponential or trigonometric forms of solutions $f(z)$ when non-trivial solutions exist.
\item[3.] Dimensional Generalization: Bridging results from a single complex variable $\mathbb{C}$ to $\mathbb{C}^m$, where spatial partial derivatives and multi-dimensional shifts $c = (c_1, \ldots, c_m)$ introduce rich geometric and analytical dynamics.
\end{enumerate}

\medskip
\noindent
We define $\mathbb{Z}_+=\mathbb{Z}[0,+\infty)=\{n\in \mathbb{Z}: 0\leq n<+\infty\}$ and $\mathbb{Z}^+=\mathbb{Z}(0,+\infty)=\{n\in \mathbb{Z}: 0<n<+\infty\}$. Also on $\mathbb{C}^m$, we define
\begin{align*}
\displaystyle \partial_{z_i}(f(z))=\frac{\partial f(z)}{\partial z_i},\ldots,\partial^{l_i}_{z_i}(f(z))=\frac{\partial^{l_i} f(z)}{\partial z_i^{l_i}},\;\;\partial^{I}(f(z))=\frac{\partial^{|I|}f(z)}{\partial z_1^{i_1}\cdots \partial z_m^{i_m}},
\end{align*}
where $l_i\in \mathbb{Z}_+$ and $i=1,2,\ldots,m$ and  $|I|=\sum_{j=1}^m i_j$.

\medskip
\noindent
In 2008, Li \cite[Theorem 1]{L1} proved that meromorphic solutions of the Fermat-type functional equation $f^2+g^2=1$ in $\mathbb{C}^2$ must be constant if and only if $\partial_{z_2}(f(z_1,z_2))$ and $\partial_{z_1}(g(z_1,z_2))$ share $0$ CM.  Moreover if $f(z_1,z_2)=\partial_{z_1}(u(z_1,z_2))$ and $g(z_1,z_2)=\partial_{z_2}(u(z_1,z_2))$, then any entire solutions of the partial differential equations:
\begin{align*}\left(\partial_{z_1}(u(z_1,z_2))\right)^2+\left(\partial_{z_2}(u(z_1,z_2))\right)^2=1\end{align*}
in $\mathbb{C}^2$ are necessarily linear (see \cite{DK1}).  The solutions of Fermat-type partial differential equations were originally investigated by Li \cite{L2} and Saleeby \cite{GS1}.

\smallskip
\noindent
Recently, many authors  (see \cite{Haldar_Mediterian_2023}, \cite{L2}, \cite{Lu2}, \cite{MSP}, \cite{GS1}, \cite{GS2}, \cite{XuHaldar2024}, \cite{XC1}) have paid attention to the existence of solutions of Fermat-type difference and partial differential-difference equations using the difference Nevanlinna theory for meromorphic functions (see \cite{TBC1, CK1, CX1, K1}) of higher dimension.

\medskip
\noindent
In the paper, we consider the following Fermat-type partial differential-difference equations
\begin{align}\label{pds1}
(P_1(z)f(z+c))^{m_1}+\left(P_2(z)f(z)+\sideset{}{_{j=1}^m}{\sum} Q_j(z)\partial_{z_j}(f(z))\right)^{m_2}=P(z),
\end{align}

\begin{align}\label{pds2} (P_1(z)f(z))^{m_1}+\left(P_2(z)f(z+c)+\sideset{}{_{j=1}^m}{\sum} Q_j(z)\partial_{z_j}(f(z))\right)^{m_2}=P(z)
\end{align}
and
\begin{align}\label{pds3} \left(\sideset{}{_{j=1}^m}{\sum} Q_j(z)\partial_{z_j}(f(z))\right)^{m_1}+\left(P_1(z)f(z+c)+P_2(z)f(z)\right)^{m_2}=P(z),
\end{align}
where $m_1,m_2\in\mathbb{N}$ such that $m_1+m_2>2$, $c=(c_1,c_2,\ldots,c_m)\in\mathbb{C}^m$, $P(\not\equiv 0)$, $P_1$, $P_2$ and $Q_i\;(i=1,2,\ldots,m)$ are polynomials in $\mathbb{C}^m$ such that $P(z)$ is irreducible in $\mathbb{C}^m$. 
\smallskip

\noindent
The purpose of this paper is to characterize all finite-order transcendental entire solutions of the equations \eqref{pds1}-\eqref{pds3} for different choices of $m_1$ and $m_2$, thereby extending Theorems A, B and C to higher dimensions.
\medskip

\noindent
The remainder of this paper is organized as follows. In Section \ref{sec2}, we recall some basic notation of Nevanlinna theory in $\mathbb{C}^m$ and collect several auxiliary lemmas for use in subsequent sections. In Section \ref{sec3}, we investigate the transcendental entire solutions of Eq. \eqref{pds1}. In Section \ref{sec4}, we determine the transcendental entire solutions of Eq. \eqref{pds2}. Furthermore, in Section \ref{sec5}, we obtain the transcendental entire solutions of Eq. \eqref{pds3}.

\section{\bf{Preliminaries on Nevanlinna theory in $\mathbb{C}^m$ and lemmas}}\label{sec2}

\noindent
Nevanlinna theory, originally developed for meromorphic functions of a single complex variable, has been extensively extended to the setting of several complex variables to investigate the value distribution of meromorphic mappings from $\mathbb{C}^m$ into complex projective spaces and complex manifolds. This theory provides deep insights into: 
(a) the growth and value distribution of meromorphic mappings; 
(b) uniqueness and rigidity phenomena for solutions of linear and nonlinear partial differential equations; and 
(c) criteria for compactness and convergence of families of meromorphic mappings. 
Its applications encompass a wide range of areas, including complex analysis, partial differential equations, complex geometry, and mathematical physics, offering both powerful analytical techniques and elegant theoretical results. The references \cite{Banerjee-Majumder-2026}-\cite{Banerjee-Majumder-Pramanik-2026}, \cite{ps1}, \cite{LB}-\cite{Majumder-Sarkar-2027}, \cite{Xu-Majumder-Sarar-Pramanik-2026} provide a comprehensive foundation for understanding the current state of research in Nevanlinna value distribution theory in several complex variables.
\smallskip

\noindent
 On $\mathbb{C}^m$, the exterior derivative $d$ splits $d= \partial+ \bar{\partial}$ and twists to $d^c= \frac{\iota}{4\pi}\left(\bar{\partial}- \partial\right)$ (see \cite{HLY1,WS1}). Clearly $dd^{c}= \frac{\iota}{2\pi}\partial\bar{\partial}$.
On $\mathbb{C}^m$, we define $\tau_m(z)=||z||^2$. The standard Kaehler metric on $\mathbb{C}^m$ is given by $\upsilon_m=dd^c\tau_m>0$. On $\mathbb{C}^m\backslash \{0\}$, we define $\omega_m=dd^c\log \tau_m\geq 0$ and $\sigma_m=d^c\log \tau_m \wedge \omega_m^{m-1}$. For any $S\subseteq \mathbb{C}^m$, let $S[r]$, $S(r)$ and $S\langle r\rangle$ be the intersection of $S$ with
respectively the closed ball, the open ball, the sphere of radius $r>0$ centered at $0\in \mathbb{C}^m$.
\smallskip

\noindent
Let $f$ be a holomorphic function on $G(\not=\varnothing)$, where $G$ is an open subset of $\mathbb{C}^m$. Then we can write $f(z)=\sum_{i=0}^{\infty}P_i(z-a)$, where the term $P_i(z-a)$ is either identically zero or a homogeneous polynomial of degree $i$. Certainly the zero multiplicity $\mu^0_f(a)$ of $f$ at a point $a\in G$ is defined by $\mu^0_f(a)=\min\{i:P_i(z-a)\not\equiv 0\}$.
\smallskip

\noindent
Let $f$ be a meromorphic function on $\mathbb{C}^m$. Then there exist holomorphic functions $g$ and $h$ such that $hf=g$ on $\mathbb{C}^m$ and $\dim_z h^{-1}(\{0\})\cap g^{-1}(\{0\})\leq m-2$. Therefore the $c$-multiplicity of $f$ is just $\mu^c_f=\mu^0_{g-ch}$ if $c\in\mathbb{C}$ and $\mu^c_f=\mu^0_h$ if $c=\infty$. The function $\mu^c_f: \mathbb{C}^m\to \mathbb{Z}$ is nonnegative and is called the $c$-divisor of $f$. 
If $f\not\equiv 0$, then $\nu=\mu_f=\mu^0_f-\mu^{\infty}_f$ is called the divisor of $f$.  For $t>0$, the counting function $n_{\nu}$ is defined by
\begin{align*}
 n_{\nu}(t)=t^{-2(m-1)}\int_{A[t]}\nu \upsilon_m^{m-1},
 \end{align*}
where $A=\text{supp}\;\nu=\ol{\{z\in G: \nu_f(z)\neq 0\}}$. The valence function of $\nu$ is defined by 
\begin{align*}
N_{\nu}(r)=N_{\nu}(r,r_0)=\int_{r_0}^r n_{\nu}(t)\frac{dt}{t}\;\;(r\geq r_0).
\end{align*}
Also we write $N_{\mu_f^a}(r)=N(r,a;f)$ if $a\in\mathbb{C}$ and $N_{\mu_f^a}(r)=N(r,f)$ if $a=\infty$.
With the help of the positive logarithm function, we define the proximity function of $f$ by
\begin{align*}
m(r, f)=\int_{\mathbb{C}^m\langle r\rangle} \log^+ |f|\;\sigma_m \geq  0.
\end{align*}
The characteristic function of $f$ is defined by $T(r, f)=m(r,f)+N(r,f)$.
We define $m(r,a;f)=m(r,f)$ if $a=\infty$ and $m(r,a;f)=m(r,1/(f-a))$ if $a\in\mathbb{C}$. Now if $a\in\mathbb{C}$, then the first main theorem becomes $m(r,a;f)+N(r,a;f)=T(r,f) + O(1)$, where $O(1)$ denotes a bounded function when $r$ is sufficiently large. We define the order of $f$ by
\begin{align*}
\rho(f):=\limsup _{r \rightarrow \infty} \frac{\log T(r, f)}{\log r}.
\end{align*}

\noindent
We define the linear measure $m(E):=\int_E dt$, the logarithmic measure $l(E):=\int_{E\cap [1,\infty)} \frac{d t}{t}$ and the upper density measure
\begin{align*}
\ol {\text{dens}}\;E=\lim\limits_{r\rightarrow \infty}\frac{1}{r}\int_{E\cap [1,r]} dt
\end{align*}
for a set $E\subset [0,\infty)$. Moreover, if $l(E)<+\infty$, resp., $m(E)<+\infty$, then $E$ is of zero upper density.
Let $S(f)=\{g:\mathbb{C}^m\to\mathbb{P}^1\;\text{meromorphic}: \parallel T(r,g)=o(T(r,f))\}$,
where $\parallel$ indicates that the relation holds only outside a set $E$ with zero upper density measure and the element in $S(f)$ is called the small function of $f$.

\noindent
To facilitate the subsequent analysis, we recall the following auxiliary lemmas.

\begin{lem}\label{L.1} \cite[Lemma 1.37]{HLY1} Let $f$ be a non-constant meromorphic function in $\mathbb{C}^m$ and $I=(\alpha_1,\ldots,\alpha_m)\in \mathbb{Z}^m_+$ be a multi-index. Then for any $\varepsilon>0$, we have
\begin{align*}
m\left(r,\partial^I(f)/f\right)\leq |I|\log^+T(r,f)+|I|(1+\varepsilon)\log^+\log T(r,f)+O(1)
\end{align*}
for all $r\not\in E$, where $l(E)<+\infty$.
\end{lem}

\begin{lem}\label{L.2} \cite[Lemma 1.2]{HY1} Let $f$ be a non-constant meromorphic function in $\mathbb{C}^m$ and let $a_1,a_2,\ldots,a_q$ be different points in $\mathbb{C}\cup \{\infty\}$. Then
\begin{align*}
 \parallel (q-2)T(r,f)\leq \sideset{}{_{j=1}^{q}}{\sum} \ol N(r,a_j;f)+O(\log (rT(r,f))).
\end{align*}
\end{lem}

\begin{lem}\label{L.3} \cite[Theorem 1.26]{HLY1} Let $f$ be non-constant meromorphic function in $\mathbb{C}^m$. Assume that 
$R(z, w)=\frac{A(z, w)}{B(z, w)}$. Then
\begin{align*}
 T\left(r, R_f\right)=\max \{p, q\} T(r, f)+O\Big(\sideset{}{_{j=0}^p}{\sum} T(r, a_j)+\sideset{}{_{j=0}^q}{\sum}T(r, b_j)\Big),
\end{align*}
where $R_f(z)=R(z, f(z))$ and two coprime polynomials $A(z, w)$ and $B(z,w)$ are given
respectively as follows: $A(z,w)=\sum_{j=0}^p a_j(z)w^j$ and $B(z,w)=\sum_{j=0}^q b_j(z)w^j$.
\end{lem}

\begin{lem}\label{L.6}\cite[Lemma 3.2]{HLY1} Let $f_j\not\equiv 0\;(j=1,2,\ldots,n;n\geq 3)$ be meromorphic functions in $\mathbb{C}^m$ such that $f_1,\ldots,f_{n-1}$ are non-constants and $f_1+\cdots+f_n=1$. If
\begin{align*}
 \parallel \;\;\sideset{}{_{j=1}^n}{\sum}\Big\lbrace N_{n-1}(r,0;f_j)+(n-1)\ol{N}(r,f_j)\Big\rbrace<\lambda T(r,f_j)+O(\log^+(T(r,f_j))
\end{align*}
holds for $j=1,2,\ldots,n-1$, where $\lambda <1$, then $f_n\equiv 1$.
\end{lem}

\begin{lem}\label{L.8}\cite[Theorem 2.1]{CX1} Let $f$ be a non-constant meromorphic function in $\mathbb{C}^m$ and let $c\in \mathbb{C}^m$. If 
\begin{align}\label{L.81}\limsup\limits_{r\rightarrow \infty} \frac{\log T(r,f)}{r}=0,
\end{align}
then
$\parallel\;m(r,f(z+c)/f(z))+m(r,f(z)/f(z+c))=o(T(r,f)).$
\end{lem}

\begin{lem}\label{L.7}\cite[Theorem 2.2]{CX1} Let $f$ be a non-constant meromorphic function in $\mathbb{C}^m$ such that \eqref{L.81} holds. 
Then $\parallel T(r,f(z+c))=T(r,f)+o(T(r,f))$ holds for any constant $c\in\mathbb{C}^m$.
\end{lem}

\noindent
Define complex differential-difference polynomials as follows:
\begin{align}\label{cl1} 
&P(f(z))\\=&\sideset{}{_{\mathbf{p} \in I}}{\sum} a_{\mathbf{p}}(z) f^{p_0}(z)\big(\partial^{\mathbf{i}_1}(f(z))\big)^{p_1} \cdots\big(\partial^{\mathbf{i}_l}(f(z))\big)^{p_l}f^{p_{l+1}}(z + \hat q_{l+1}) \cdots f^{p_{l+s}}(z + \hat q_{l+s}),\nonumber
\end{align}
$\mathbf{p}=\left(p_0, \ldots, p_{l+s}\right) \in \mathbb{Z}_{+}^{l+s+1}$,
\begin{align}\label{cl2} 
&Q(f(z))\\=&\sideset{}{_{\mathbf{q} \in J}}{\sum}c_{\mathbf{q}}(z) f^{q_0}(z)\big(\partial^{\mathbf{j}_1}(f(z))\big)^{q_1} \cdots\big(\partial^{\mathbf{j}_l}(f(z))\big)^{q_l}f^{q_{l+1}}(z+\tilde q_{l+1}) \cdots f^{q_{l+t}}(z+\tilde q_{l+t}),\nonumber
\end{align}
$\mathbf{q}=\left(q_0, \ldots, q_{l+t}\right) \in \mathbb{Z}_{+}^{l+t+1}$ and 
\begin{align}\label{cl3} 
B(f(z))=\sideset{}{_{k=0}^{n}}{\sum} b_k(z) f^k(z),
\end{align}
where $I, J$ are finite sets, $a_{\mathbf{p}}, b_k, c_{\mathbf{q}}\in S(f)$ such that $b_n\not\equiv 0$ and $\hat q_i, \tilde q_j\in\mathbb{C}^m$. 
\smallskip

\noindent
Now by using Lemmas \ref{L.1} and \ref{L.8} and proceeding in the same way as done in the proofs of Lemma 2.1 \cite{ps1} and Theorem 2.3 \cite{LY1}, we get the following lemma.

\begin{lem}\label{L.4} Let $f(z)$ be a non-constant meromorphic function in $\mathbb{C}^m$ such that \eqref{L.81} holds.
Assume that $f$ satisfies the differential-difference equation $B(f)Q(f)=P(f)$,
where $P(f)$, $Q(f)$ and $B(f)$ are defined as in (\ref{cl1}), (\ref{cl2}) and (\ref{cl3}) respectively.
If $\deg(P(f)) \leq n=\deg(B(f))$, then $\parallel m(r, Q(f))=o(T(r, f)).$
\end{lem}

\begin{lem}\label{L.4e}\cite[Proposition 3.2]{hy1} Let $P$ be a non-constant entire function in $\mathbb{C}^n$. Then
	\[\rho(e^P)=
	\begin{cases}
		\deg(P) & \text{if $P$ is a polynomial,}\\
		+\infty & \text{otherwise}.
	\end{cases}\]
\end{lem}

\noindent
The following lemma follows directly from the proof of Theorem 1 \cite{Y1}.
\begin{lem}\label{L.5}
Let $f$ and $g$ be two non-constant entire functions in $\mathbb{C}^m$ and let $m_1$ and $m_2$ be two positive integers. Then the functional equation $f^{m_1}(z)+g^{m_2}(z)=1$
holds only when $1/m_1+1/m_2\geq 1$.
\end{lem}

\section {{\bf The transcendental entire solutions of Eq. (\ref{pds1})}}\label{sec3} 

\noindent
In this section, we present our main result, which provides a multidimensional analogue of Theorem A.

\begin{theo}\label{t2.1} Let $P(\not\equiv 0)$, $P_1$, $P_2$ and $Q_i\;(i=1,2,\ldots,m)$ are polynomials in $\mathbb{C}^m$ such that $P_1Q_j\not\equiv 0\;(j=1,2,\ldots,m)$ and $P(z)$ is irreducible in $\mathbb{C}^m$.  Then any finite order transcendental entire solutions of Eq. (\ref{pds1}) are characterized as follows:
	\begin{enumerate}
		\item[\em{(i)}] if $m_1=m_2$, then $m_1=m_2=2$. Also
		\begin{align*}\displaystyle
			f(z)=k\frac{\exp(\iota \ul{P^{*}}(z))-\ul{P}(z)\exp(-\iota \ul{P^{*}}(z))}{2\iota \ul{P_1}(z)},
		\end{align*}
		where $P^*(z)=d_1z_1+d_2z_2+\ldots+d_mz_m+\tilde P^*(z)$, $d_1,d_2,\ldots,d_m\in\mathbb{C}$, $\tilde P^*(z)$ is a polynomial in $\mathbb{C}^m$ such that $\tilde P^*(z-c)\equiv \tilde P^*(z)$ and
		\begin{align*}\displaystyle
			e^{-\iota(c_1d_1+\ldots+c_md_m)} \left(-\iota P_2\ul{P_1}+k\sideset{}{_{j=1}^{m}}{\sum}Q_j\left(\ul{P_1}\partial_{z_j}(\ul{P^*})+\iota\partial_{z_j}(\ul{P_1})\right)\right)\equiv \ul{P_1}^2,\;(k=\pm 1);
		\end{align*} 
	\smallskip

\item[\em{(ii)}] if $m_2>m_1$, $\sum_{j=1}^m \partial_{z_j}(P_1(z))=P_2(z+c)$ and $Q_j(z+c)=P_1(z)$ for $j=1,2,\ldots,m$, then $m_1=1$, $m_2=2$ and 
\begin{align*}\displaystyle
f(z)=\frac{1}{\ul{P_1}(z)}\left(-\left(-\sideset{}{_{j=1}^m}{\sum}(z_j-c_j)/2m+\ul{g_1}(z)\right)^2+\ul{P}(z)\right),
\end{align*}
where $g_1(z)=g_1(z_2-z_1,z_3-z_2,\ldots,z_m-z_{m-1})$ is an arbitrary finite-order transcendental entire function such that
\begin{align*} 
\displaystyle
g_1(z+c)-g_1(z)=(c_1+c_2+\ldots+c_m)/2m+\sideset{}{_{j=1}^m}{\sum}\partial_{z_j}(P(z));
\end{align*}

\item[\em{(iii)}] if $m_1>m_2$, then Eq. (\ref{pds1}) does not have any finite-order transcendental entire solution.
\end{enumerate}
\end{theo}

As a special case, choosing $P_1 =P = 1$, $p_2 = 0$, $Q_j=1\;(j=1,2,\ldots,m)$ and $m_1=m_2=2$ in Theorem \ref{t2.1} leads to the following corollary, which improves \cite[Theorem 1.1]{XC2}.
\begin{cor} All finite-order transcendental entire solutions of 
\begin{align*}
\left(\sideset{}{_{j=1}^m}{\sum}\partial_{z_j}(f(z))\right)^2+f^2(z+c)=1
\end{align*}
must satisfy
\begin{align*}\displaystyle
f(z)=\frac{\exp(\iota \ul{P^{*}}(z))-\exp(-\iota \ul{P^{*}}(z))}{2\iota},
\end{align*}
where $P^*(z)=d_1z_1+d_2z_2+\ldots+d_mz_m+\tilde P^*(z)$, $d_1,d_2,\ldots,d_m\in\mathbb{C}$, $\tilde P^*(z)$ is a polynomial in $\mathbb{C}^m$ such that $\tilde P^*(z-c)\equiv \tilde P^*(z)$ and
	\begin{align*}\displaystyle
		\exp\left(-\iota(c_1d_1+\ldots+c_md_m)\right)\sideset{}{_{j=1}^{m}}{\sum}\partial_{z_j}(\ul{P^*})=1.
	\end{align*} 
\end{cor}

\begin{proof}[{\bf Proof of Theorem \ref{t2.1}}]
Let $f(z)$ be any finite-order transcendental entire solution of Eq. (\ref{pds1}). Set 
\begin{align}\label{sm}
\displaystyle
F_f(z)=P_2(z)f(z)+\sideset{}{_{j=1}^m}{\sum} Q_j(z)\partial_{z_j}(f(z)).
\end{align}

\noindent
If $F_f(z)$ is a constant, then \eqref{pds1} implies that $P_1(z)f(z)$ is also a constant, which is impossible since $f(z)$ is transcendental. Hence $F_f(z)$ is non-constant. Using Lemma \ref{L.3} to (\ref{pds1}), we get
\begin{align}\label{sm.0} \parallel\;m_2T(r,F_f)+o(T(r,F_f))=m_1T(r,f(z+c))+o(T(r,f(z+c))).
\end{align}
Now using Lemma \ref{L.5} to (\ref{pds1}), we get
\begin{align}\label{sm.1} 
1/m_1+1/m_2\geq 1.
\end{align}

\noindent
We consider the following cases.\par
\smallskip

\noindent
\noindent
{\bf Case 1.} Let $m_1=m_2$. Then (\ref{sm.1}) gives $m_1=m_2=2$ and so $F_f^2(z)+(P_1(z)f(z+c))^2=P(z)$, i.e.,
\begin{align}\label{sm.2} \displaystyle
	(F_f(z)+\iota P_1(z) f(z+c))(F_f(z)-\iota P_1(z) f(z+c))=P(z).
\end{align}
Therefore we may assume that 
\begin{align}\label{sm.3}\displaystyle
	\begin{cases}
		F_f(z)+\iota P_1(z)f(z+c)=\exp(\iota P^{*}(z))\\
		F_f(z)-\iota P_1(z)f(z+c)=P(z)\exp(-\iota P^{*}(z)),
	\end{cases}
\end{align}
or 
\begin{align}\label{sm.4}\displaystyle
	\begin{cases}
		F_f(z)-\iota P_1(z)f(z+c)=\exp(\iota P^{*}(z))\\
		F_f(z)+\iota P_1(z)f(z+c)=P(z)\exp(-\iota P^{*}(z)),\end{cases}
\end{align}
where $P^{*}(z)$ is an entire function in $\mathbb{C}^m$. Solving, we get
\begin{align}\label{sm.5}\displaystyle 
	F_f(z)=P_2(z)f(z)+\sideset{}{_{j=1}^m}{\sum} Q_j(z)\frac{\partial f(z)}{\partial z_j}=\frac{\exp(\iota P^{*}(z))+P(z)\exp(-\iota P^{*}(z))}{2}
\end{align}
and 
\begin{align*} \displaystyle
	f(z+c)=k\frac{\exp(\iota P^{*}(z))-P(z)\exp(-\iota P^{*}(z))}{2\iota P_1(z)},
\end{align*}
i.e., 
\begin{align}\label{sm.6} \displaystyle
	f(z)=k\frac{\exp(\iota \ul{P^{*}}(z))-\ul{P}(z)\exp(-\iota \ul{P^{*}}(z))}{2\iota \ul{P_1}(z)},
\end{align}
where $k=1$ for \eqref{sm.3} and $k=-1$ for \eqref{sm.4}.
Now using Lemma \ref{L.3} to \eqref{sm.6}, we get $\rho\left(\exp(\iota P^{*}(z-c))\right)<+\infty$. Consequently by Lemma \ref{L.4e}, we conclude that $P^{*}(z)$ is a non-constant polynomial in $\mathbb{C}^m$. Again from (\ref{sm.6}), we have
\begin{align*} \displaystyle
	\partial_{z_j}(f(z))=&k\left(\frac{1}{2\ul{P_1}(z)}\partial_{z_j}(\ul{P^*}(z))-\frac{1}{2 \iota \ul{P_1}^2(z)}\partial_{z_j}(\ul{P_1}(z))\right)\exp (\iota \ul{P^{*}}(z))\\+&k\left(\frac{\ul{P}(z)}{2\ul{P_1}(z)}\partial_{z_j}(\ul{P^*}(z))+\frac{\ul{P}(z)}{2\iota \ul{P_1}^2(z)}\partial_{z_j}(\ul{P_1}(z))-\frac{1}{2 \iota \ul{P_1}(z)}\partial_{z_j}(\ul{P}(z))\right) \exp(-\iota \ul{P^{*}}(z)),
\end{align*}
where $j=1,2,\ldots,m$ and so 
\begin{align}\label{sm.7}\displaystyle
	F_f(z)=&P_2(z)f(z)+\sideset{}{_{j=1}^m}{\sum} Q_j(z)\partial_{z_j}(f(z))\\=&\nonumber\left(\frac{P_2}{2\iota \ul{P_1}}+k\sideset{}{_{j=1}^{m}}{\sum}Q_j\left(\frac{1}{2\ul{P_1}}\partial_{z_j}(\ul{P^*})-\frac{1}{2 \iota \ul{P_1}^2}\partial_{z_j}(\ul{P_1})\right)\right)\exp(\iota \ul{P^{*}})\\+&\nonumber\left(-\frac{P_2\ul{P}}{2 \iota \ul{P_1}}+k\sideset{}{_{j=1}^{m}}{\sum}Q_j\left(\frac{\ul{P}}{2\ul{P_1}}\partial_{z_j}(\ul{P^*})+\frac{\ul{P}}{2\iota \ul{P_1}^2}\partial_{z_j}(\ul{P_1})-\frac{1}{2 \iota \ul{P_1}}\partial_{z_j}(\ul{P})\right)\right)\exp(-\iota \ul{P^{*}}).\nonumber
\end{align}
Clearly from (\ref{sm.5}) and (\ref{sm.7}), we get
\begin{align}\label{sm17} 
	\alpha \exp(\iota \ul{P^{*}}(z))+\beta \exp(-\iota \ul{P^{*}}(z))-\exp(\iota P^{*}(z))-P\exp(-\iota P^{*}(z))\equiv 0, 
\end{align} 
where 
\begin{align}\label{sm17a}\displaystyle
	\alpha\equiv  \left(\frac{P_2}{\iota \ul{P_1}}+k\sideset{}{_{j=1}^{m}}{\sum}Q_j\left(\frac{1}{\ul{P_1}}\partial_{z_j}(\ul{P^*})-\frac{1}{ \iota \ul{P_1}^2}\partial_{z_j}(\ul{P_1})\right)\right)
\end{align} 
and 
\begin{align}\label{sm17b}\displaystyle
	\beta \equiv\left(-\frac{P_2\ul{P}}{ \iota \ul{P_1}}+k\sideset{}{_{j=1}^{m}}{\sum}Q_j\left(\frac{\ul{P}}{\ul{P_1}}\partial_{z_j}(\ul{P^*})+\frac{\ul{P}}{\iota \ul{P_1}^2}\partial_{z_j}(\ul{P_1})-\frac{1}{ \iota \ul{P_1}}\partial_{z_j}(\ul{P})\right)\right). 
\end{align}
Now, from (\ref{sm17}), we find
\begin{align}\label{sm18}
	\alpha\exp (\iota (\ul{P^*}-P^*))+\beta\exp (-\iota (\ul{P^*}+P^*)) -P\exp  (-2\iota P^*)\equiv 1.
\end{align}
Since $P^*(z)$ is a non-constant polynomial in $\mathbb{C}^m$, it is easy to verify that $\ul{P^*}+P^*$ is also a non-constant polynomial in $\mathbb{C}^m$. 
Now we consider following four sub-cases.\par
\medskip

\noindent
{\bf Sub-case 1.1.} Assume that $\alpha\equiv \beta \equiv 0.$ Then (\ref{sm18}) reduces to $P\exp  (-2\iota P^*)\equiv -1$, which is impossible by Lemma \ref{L.3}.
\medskip

\noindent
{\bf Sub-case 1.2.}
Assume that $\alpha\equiv0$ and $\beta \not\equiv 0.$ Then (\ref{sm18}) reduces to 
\begin{align}\label{sm19}
	\beta \exp (-\iota (\ul{P^*}+P^*)) -P\exp  (-2\iota P^*)\equiv 1.
\end{align} 
Clearly $T(r,\beta)=o(r, \exp  (-2\iota P^*))$, $T(r,\beta)=o(r, \exp (-\iota (\ul{P^*}+P^*)))$, $T(r,P)=o(r,\exp P^*)$ and $T(r,P)=o(r, \exp (-\iota (\ul{P^*}+P^*)))$. 

\noindent
Let $F=\beta \exp (-\iota (\ul{P^*}+P^*))$ and $G=P\exp  (-2\iota P^*)$. Using Lemma \ref{L.2}, we get
\begin{align*}\displaystyle
	T(r,F) \le& \ol N(r, F)+\ol N\left(r,0;F\right)+\ol N\left(r,1;F\right) +o(T(r,F))\\=& \ol N(r, F)+\ol N\left(r,0;F\right)+\ol N\left(r,0;G\right) +o(T(r,F))\\\le & \ol N(r, \beta)+\ol N\left(r,0;\beta\right)+\ol N\left(r,0;P\right) +o(T(r,F)) \\\le&  2T(r, \beta)+T\left(r,P\right) +o(T(r,F))=o(T(r,F)),
\end{align*}
which leads to a contradiction.
\medskip

\noindent
{\bf Sub-case 1.3.}
Assume that $\alpha\not\equiv0$ and $\beta \equiv 0.$ Then (\ref{sm18}) reduces to
\begin{align}\label{sm19a}
	\alpha \exp (\iota (\ul{P^*}-P^*))-P\exp  (-2\iota P^*)\equiv 1.
\end{align}
If $\ul{P^*}-P^*$ is constant, then using Lemma \ref{L.3} to \eqref{sm19a}, we get a contradiction. Hence $\ul{P^*}-P^*$ is non-constant.  Now, proceeding in the same way as done in the proof of Sub-case 1.2, we get a contradiction.
\medskip

\noindent
{\bf Sub-case 1.4.}
Assume that $\alpha\not\equiv0$ and $\beta \not\equiv 0.$ Then applying Lemma \ref{L.6} to (\ref{sm18}), we get 
\begin{align}\label{sm19b}
	\alpha \exp (\iota (\ul{P^*}-P^*))\equiv 1,
\end{align}
which shows that $\ul{P^*}-P^*$ is a constant and so we may assume that
\begin{align}\label{sm19c}
	P^*(z)=d_1z_1+d_2z_2+\ldots+d_mz_m+\tilde P^*(z),
\end{align} 
where $d_1,d_2,\ldots,d_m\in\mathbb{C}$ and $\tilde P^*(z)$ is a polynomial in $\mathbb{C}^m$ such that $\tilde P^*(z-c)\equiv \tilde P^*(z)$.
Substituting \eqref{sm19b} into \eqref{sm18}, we get
\begin{align}\label{sm19d}
	\beta\exp (-\iota (\ul{P^*}-P^*))\equiv P.
\end{align}
Let $d=\exp (\iota (\ul{P^*}-P^*))$. Now from \eqref{sm19b} and \eqref{sm19d}, we get $d\alpha=1$, $\beta=dP$ and
$\alpha\beta \equiv P$. Therefore from \eqref{sm17a} and \eqref{sm17b}, we obtain
\begin{align}\label{sm19e}\displaystyle
	d \left(-\iota P_2\ul{P_1}+k\sideset{}{_{j=1}^{m}}{\sum}Q_j\left(\ul{P_1}\partial_{z_j}(\ul{P^*})+\iota\partial_{z_j}(\ul{P_1})\right)\right)\equiv \ul{P_1}^2.
\end{align} 
%and 
%\begin{align}\label{sm19g}\displaystyle
%&\left(-P_2\ul{P_1}+\sum_{j=1}^{m}Q_j\left(\ul{P_1}\frac{\partial \ul{P^*}}{\partial z_j}+\frac{\partial \ul{P_1}}{\partial z_j}\right)\right)
%\left(P_2\ul{P} \ul{P_1}+\sum_{j=1}^{m}Q_j\left(\ul{P}\ul{P_1}\frac{\partial \ul{P^*}}{\partial z_j}-\ul{P}\frac{\partial \ul{P_1}}{\partial %z_j}+\ul{P_1}\frac{\partial \ul{P}}{\partial z_j}\right)\right)\nonumber\\\equiv& P\ul{P_1}^4.
%\end{align}
Finally from \eqref{sm.6}, we have
\begin{align*}\displaystyle
	f(z)=k\frac{\exp(\iota \ul{P^{*}}(z))-\ul{P}(z)\exp(-\iota \ul{P^{*}}(z))}{2\iota \ul{P_1}(z)},
\end{align*}
where $P^*(z)=d_1z_1+d_2z_2+\ldots+d_mz_m+\tilde P^*(z)$, $d_1,d_2,\ldots,d_m\in\mathbb{C}$, $\tilde P^*(z)$ is a polynomial in $\mathbb{C}^m$ such that $\tilde P^*(z-c)\equiv \tilde P^*(z)$ and
\begin{align*}\displaystyle
	\exp\left(-\iota(c_1d_1+\ldots+c_md_m)\right) \left(-\iota P_2\ul{P_1}+k\sideset{}{_{j=1}^{m}}{\sum}Q_j\left(\ul{P_1}\partial_{z_j}(\ul{P^*})+\iota\partial_{z_j}( \ul{P_1})\right)\right)\equiv \ul{P_1}^2.
\end{align*} 
%and \eqref{sm19g} holds.
\smallskip

\noindent
{\bf Case 2.} Let $m_1\neq m_2$. We now consider following two sub-cases.
\smallskip

\noindent
{\bf Sub-case 2.1.} Let $m_1>m_2$. If $m_2\geq 2$, then from (\ref{sm.1}), we get a contradiction. Hence $m_2=1$ and so from (\ref{pds1}), we get
\begin{align}\label{al1} 
	\left(P_1(z)f(z+c)\right)^{m_1}+F_f(z)=P(z).
\end{align}
Since $\rho(f)<+\infty$, \eqref{L.81} holds 
and so by Lemma \ref{L.7}, we obtain
\begin{align}\label{sm.11} 
	\parallel\;T(r,f(z+c))=T(r,f)+o(T(r,f)).
\end{align}
Now using Lemmas \ref{L.1}, \ref{L.3} and (\ref{sm.11}) to (\ref{al1}), we get
\begin{align*}\displaystyle
	\parallel\;m_1 T(r, f))=& m_1 T(r, P_1(z)f(z+c))+o(T(r, f)) \nonumber\\
	=&T\left(r, \left(P_1(z)f(z+c)\right)^{m_1}\right)+o(T(r, f))\nonumber\\
	=&T\left(r,F_f(z)-P(z)\right)+o(T(r, f)) \nonumber\\
	%=&T\left(r,F_f(z)\right)+o(T(r, f)) \nonumber\\
	\leq&m\left(r,F_f(z)\right)+o(T(r, f)) \nonumber\\
	\leq &m\left(r, \sideset{}{_{j=1}^m}{\sum} Q_j(z)\partial_{z_j}(f(z))/f(z)\right)+m(r, f(z))+o(T(r, f))\nonumber \\
	=&T(r, f)+o(T(r, f))
\end{align*}
and so $\parallel (m_1-1) T(r, f) \leq o(T(r, f)).$
Since $m_1\geq 2$, we get a contradiction.

\smallskip

\noindent
{\bf Sub-case 2.2.} Let $m_2>m_1$, $\sum_{j=1}^m \partial_{z_j}(P_1(z))=\ol{P_2}(z)$ and $\ol{Q_j}(z)=P_1(z)$ for $j=1,2,\ldots,m$.  If $m_1\geq 2$, then from (\ref{sm.1}), we get a contradiction. Hence $m_1=1$ and so from (\ref{pds1}), we get
\begin{align}\label{sm.12} 
P_1(z)f(z+c)+F_f^{m_2}(z)=P(z),
\end{align}
Differentiating (\ref{sm.12}) partially with respect to $z_j$, we get
\begin{align*} \displaystyle
m_2F_f^{m_2-1}(z)\partial_{z_j}(F_f(z))+\partial_{z_j}(P_1(z))\ol{f}(z)+P_1(z)\partial_{z_j}(\ol{f}(z))=\partial_{z_j}(P(z)),
\end{align*}
i.e.,
\begin{align}\label{sm.13} \displaystyle
&m_2F_f^{m_2-1}\sideset{}{_{j=1}^m}{\sum} \partial_{z_j}(F_f(z))+
\sideset{}{_{j=1}^m}{\sum} \partial_{z_j}(P_1(z)) \ol{f}(z)+\sideset{}{_{j=1}^m}{\sum} P_1(z)\partial_{z_j}(\ol{f}(z))
\\=&\sideset{}{_{j=1}^m}{\sum} \partial_{z_j}(P(z)).\nonumber
\end{align}
Now using (\ref{sm}) to \eqref{sm.13}, we get
\begin{align}\label{sm.14}  
m_2F_f^{m_2-1}(z)\tilde \beta(z)+\ol{F_f}(z)=\tilde \beta_1(z),
\end{align}
where
\begin{align}\label{sm.14a}
\tilde \beta(z)=\sideset{}{_{j=1}^m}{\sum} \partial_{z_j}(F_f(z))
\end{align}
and 
\begin{align}\label{sm.14c}
\tilde \beta_1(z)=\sideset{}{_{j=1}^m}{\sum} \partial_{z_j}(P(z)).
\end{align}
If possible, suppose $\tilde \beta\equiv 0$. Then from (\ref{sm.14}), we have 
\begin{align*}\displaystyle
\ol{F_f}(z)=\sideset{}{_{j=1}^m}{\sum}\partial_{z_j}(P(z)),
\end{align*}
which shows that $F_f(z)$ is a polynomial in $\mathbb{C}^m$. Therefore it is easy to conclude from \eqref{sm.12} that $f(z)$ is a polynomial in $\mathbb{C}^m$, which contradicts the fact that $f(z)$ is transcendental. 
Hence $\tilde \beta\not\equiv 0$.
We know that $m_2\geq 2$. Now using Lemma \ref{L.4} to (\ref{sm.14}), we get
$\parallel m (r,\tilde \beta)=o(T(r, F_f)).$
Since $\beta$ is an entire function, we have 
\begin{align}\label{sm.15a}
\parallel T(r,\tilde \beta)=o(T(r, F_f)).
\end{align}

\noindent
Based on the above estimations, we now consider the following sub-cases:
\smallskip

\noindent
{\bf Sub-case 2.2.1.} Let $m_2-1\geq 2$. Now using Lemma \ref{L.4} to (\ref{sm.14}), we get $\parallel m(r,\tilde \beta F_f)=o(T(r,F_f))$.
Since $F_f$ is an entire function and $\parallel T(r,\tilde \beta)=o(T(r,F_f))$, we have 
\begin{align*}
\parallel\;T(r,F_f)=m(r,F_f)\leq m(r,\tilde \beta F_f)+m(r,1/\tilde \beta)\leq o(T(r,F_f))
\end{align*}
and so from (\ref{sm.0}), we have $\parallel T(r,\ol{f}(z))=o(T(r,\ol{f}(z)))$. Finally from (\ref{sm.11}), we get a contradiction.\par
\smallskip

\noindent
{\bf Sub-case 2.2.2.} Let $m_2-1=1$. Then from (\ref{sm.14}), we get
\begin{align}\label{sm.15} 
2F_f(z)\tilde \beta(z)+\ol{F_f}(z)=\beta_1(z).
\end{align}
Differentiating (\ref{sm.15}) partially with respect to $z_j$, we get
\begin{align*}
 2\partial_{z_j}(F_f(z))\tilde \beta(z)+2F_f(z)\partial_{z_j}(\tilde \beta(z))+\partial_{z_j}(\ol{F_f}(z))=0,
\end{align*}
i.e.,
\begin{align}\label{sm.17}
2\tilde \beta(z)\sideset{}{_{j=1}^m}{\sum} \partial_{z_j}(F_f(z))+2F_f(z)\sideset{}{_{j=1}^m}{\sum}\partial_{z_j}(\tilde \beta(z))+\sideset{}{_{j=1}^m}{\sum} \partial_{z_j}(\ol{F_f}(z))=0.
\end{align}
Now using (\ref{sm}) and \eqref{sm.14a} to \eqref{sm.17}, we get
\begin{align}\label{sm.18} 
2\tilde \beta^2(z)+2F_f(z)\sideset{}{_{j=1}^m}{\sum}\partial_{z_j}(\tilde \beta(z))+\ol {\tilde{\beta}}(z)=0.
\end{align}
\smallskip

\noindent
Suppose $\sum_{j=1}^m\partial_{z_j}(\tilde \beta(z))\not\equiv 0$. In view of (\ref{sm.15a}) and using Lemma \ref{L.3} to (\ref{sm.18}), we get  $\parallel T(r, F_f)=o(T(r, F_f))$ and so from (\ref{sm.0}), we get 
$\parallel T(r,\ol{f}(z))=o(T(r,\ol{f}(z)))$. Finally from (\ref{sm.11}), we get a contradiction.
Hence $\sum_{j=1}^m\partial_{z_j}(\tilde \beta(z))\equiv 0$. Then from (\ref{sm.18}), we have
\begin{align}\label{sm.20} 
2\tilde \beta^2(z)=-\ol{\tilde{\beta}}(z).
\end{align}

\noindent
If $\tilde \beta(z)$ is a transcendental entire function, then using Lemma \ref{L.4} to (\ref{sm.20}), we get 
$\parallel m(r,\tilde \beta)=o(T(r, \tilde \beta))$. This implies $\parallel T(r,\tilde \beta)=o(T(r,\tilde \beta))$, which is impossible. 
Hence $\tilde \beta(z)$ is a polynomial in $\mathbb{C}^m$.  It is easy to deduce from (\ref{sm.20}) that $\tilde \beta(z)$ is a constant, say $k_1(\neq 0)$. Then from (\ref{sm.20}), we have 
$2k_1^2=-k_1$, which implies that $k_1=-1/2$. Therefore 
\begin{align}\label{sm.21} 
\tilde \beta(z)=\sideset{}{_{j=1}^m}{\sum}\partial_{z_j}(F_f(z))=-1/2.
\end{align} 
Again from \eqref{sm.14c} and (\ref{sm.15}), we get 
\begin{align}\label{Sm.21}
F_f(z)=\ol{F_f}(z)-\sideset{}{_{j=1}^m}{\sum}\partial_{z_j}(P(z)).
\end{align}
The Lagrange auxiliary equations corresponding to (\ref{sm.21}) are as follows
\begin{align}\label{sm.21a} 
	\displaystyle
	\frac{d z_1}{1}=\frac{dz_2}{1}=\ldots=\frac{dz_m}{1}=\frac{dF_f(z)}{-1/2}=\frac{dz_1+dz_2+\ldots+dz_m}{m}.
\end{align}
Clearly $z_2-z_1=d_1$, $z_3-z_2=d_2$, $\ldots$, $z_m-z_{m-1}=d_{m-1}$ and $mF_f(z)+\frac{1}{2}(z_1+z_2+\ldots+z_m)=d_m$ be independent solutions of (\ref{sm.21a}). Therefore the general solution of (\ref{sm.21}) is as follows
\begin{align}\label{sm.22} 
\displaystyle
F_f(z)=-(z_1+z_2+\ldots+z_m)/2m+g_1(z),
\end{align}
where $g_1(z)=g_1(z_2-z_1,z_3-z_2,\ldots,z_m-z_{m-1})$ is an arbitrary finite-order transcendental entire function. Using \eqref{Sm.21} to (\ref{sm.22}), we have
\begin{align*}
	\displaystyle
	&-(z_1 + c_1+z_2+c_2+\ldots+z_m+c_m)/2m + g_1(z+c)-\sideset{}{_{j=1}^m}{\sum}\partial_{z_j}(P(z))\\=&-(z_1+z_2+\ldots+z_m)/2m + g_1(z),
\end{align*}
which implies
\begin{align}\label{sm.23} 
	\displaystyle
	g_1(z+c)-g_1(z)=(c_1+c_2+\ldots+c_m)/2m+\sideset{}{_{j=1}^m}{\sum}\partial_{z_j}(P(z)).
\end{align}
Finally from \eqref{sm.12}, we have 
\begin{align*}\displaystyle
f(z+c)=\frac{1}{P_1(z)}\left(-\Big(-(z_1+z_2+\ldots+z_m)/2m+g_1(z)\Big)^2+P(z)\right),
\end{align*}
where $g_1(z)=g_1(z_2-z_1,z_3-z_2,\ldots,z_m-z_{m-1})$ is an arbitrary finite-order transcendental entire function such that (\ref{sm.23}) holds.
\end{proof}

\section {{\bf The transcendental entire solutions of Eq. \eqref{pds2}}}\label{sec4}

\noindent
In this section, we state our main result, extending Theorem B to higher dimensions.

\begin{theo}\label{t2.2} Let $P(\not\equiv 0)$, $P_1$, $P_2$ and $Q_i\;(i=1,2,\ldots,m)$ are polynomials in $\mathbb{C}^m$ such that $P_1P_2\not\equiv 0$ and  $P(z)$ is irreducible in $\mathbb{C}^m$. Then any finite-order transcendental entire solutions of Eq. (\ref{pds2}) are characterized as follows:
\begin{enumerate}
\item[\em{(i)}] if $m_1=m_2$, then $m_1=m_2=2$. Also
\begin{align*}
	f(z)=k\frac{\exp(\iota Q^{*}(z))-P(z)\exp(-\iota Q^{*}(z))}{2\iota P_1(z)},
\end{align*}
where $P^*(z)=d_1z_1+d_2z_2+\ldots+d_mz_m+\tilde P^*(z)$, $d_1,d_2,\ldots,d_m\in\mathbb{C}$, $\tilde P^*(z)$ is a polynomial in $\mathbb{C}^m$ such that $\tilde P^*(z+c)\equiv \tilde P^*(z)$ and
\begin{align*}
	P_2e^{-\iota (d_1c_1+d_2c_2+\ldots+d_mc_m)}+\ol {P_1} \sideset{}{_{j=1}^m}{\sum} Q_j\left(\frac{\iota}{P_1} \partial_{z_j}(Q^*)-\frac{1}{P_1^2}\partial_{z_j}( P_1)\right)\equiv\frac{\iota  \ol{P_1}}{k},\;(k=\pm 1);
\end{align*} 
%and 
%\begin{align*}
%&\left(\sum_{j=1}^{m}Q_j\left(\frac{{P}}{2P_1}\frac{\partial {P^*}}{\partial z_j}+\frac{{P}}{2\iota {P_1}^2}\frac{\partial {P_1}}{\partial z_j}-\frac{1}{2 \iota %{P_1}}\frac{\partial {P}}{\partial z_j}\right)-\frac{kP}{2}\right)\times\\&\left(\sum_{j=1}^{m}Q_j\left(\frac{1}{2{P_1}}\frac{\partial {P^*}}{\partial z_j}-\frac{1}{2 %\iota {P_1}^2}\frac{\partial {P_1}}{\partial z_j}\right)-\frac{k}{2}\right) \equiv\frac{P_2^2\ol P}{4\ol{P_1}^2}.
%\end{align*}
\smallskip 
		
\item[\em{(ii)}] if $m_2>m_1$, $\sum_{j=1}^m \partial_{z_j}(P_1(z))=P_2(z)$, $Q_i(z)=P_1(z)$ for $j=1,2,\ldots,m$ and $\sum_{j=1}^m Q_j(z)\partial_{z_j}(f(z))$ is $c$-periodic, then $m_1=1$, $m_2=2$ and
\begin{align*}
f(z)=\frac{1}{P_1(z)}\left(-\left(-\sideset{}{_{j=1}^m}{\sum}z_j/2m+\ul{g_1}(z)\right)^2+P(z)\right),
\end{align*}
where $g_1(z)=g_1(z_2-z_1,z_3-z_2,\ldots,z_m-z_{m-1})$ is an arbitrary finite-order transcendental entire function such that
\begin{align*} 
\displaystyle
g_1(z+c)-g_1(z)=(c_1+c_2+\ldots+c_m)/2m-\sideset{}{_{j=1}^m}{\sum}\partial_{z_j}(\ol{P}(z)).
\end{align*}
\item[\em{(iii)}] if $m_1>m_2$, then Eq. (\ref{pds2}) does not have any finite-order transcendental entire solution.
\end{enumerate}
\end{theo}

\noindent
Setting $P_1 = P_2 = P = 1$, $Q_j=0\;(j=1,2,\ldots,m)$ and $m_1=m_2=2$ in Theorem \ref{t2.2} yields the following corollary, which improves upon \cite[Theorem 1.4]{XC1}.
\begin{cor} All finite-order transcendental entire solutions of $f^2(z)+f^2(z+c)=1$ must satisfy
\begin{align*}
	f(z)=\frac{\exp(\iota P^{*}(z))-\exp(-\iota P^{*}(z))}{2\iota},
\end{align*}
where $P^*(z)=d_1z_1+d_2z_2+\ldots+d_mz_m+\tilde P^*(z)$, $d_1,d_2,\ldots,d_m\in\mathbb{C}$, $\tilde P^*(z)$ is a polynomial in $\mathbb{C}^m$ such that $\tilde P^*(z+c)\equiv \tilde P^*(z)$ and $\exp({-\iota (d_1c_1+d_2c_2+\ldots+d_mc_m)})=\iota.$
\end{cor}

\noindent
The following example illustrates our findings
\begin{exm} Let 
\begin{align*}
f(z) =\frac{\exp(\iota P(z))-\exp(-\iota P(z))}{2\iota},
\end{align*}
where $P(z)=z_1+z_2+z_3+(z_1+z_2-z_3)^2+(z_1-z_2)^2$ and $c=(c_1,c_2,c_3)=(\frac{\pi }{8},\frac{\pi }{8},\frac{\pi }{4})$. Note that $P(z+c)-P(z)=\frac{\pi }{2}$ and  $f(z_1, z_2, z_3)$ is a finite-order transcendental entire function. It is easy to verify that $f(z_1, z_2, z_3)$ is a solution to the Eq. $f^2(z)+f^2(z+c)=1$.
\end{exm}

\begin{proof}[{\bf Proof of Theorem \ref{t2.2}}] Let $f(z)$ be any finite-order transcendental entire solution to the Eq. (\ref{pds2}). Take
\begin{align}\label{sm.2.1}
G_f(z)=P_2(z)f(z+c)+\sideset{}{_{j=1}^m}{\sum} Q_j(z)\partial_{z_j}(f(z)).
\end{align}

\noindent	
It is easy to prove that $G_f(z)$ is non-constant. Now, using Lemma \ref{L.3} to (\ref{pds2}), we get
\begin{align}\label{sm.2.2} 
\parallel\;m_2T(r,G_f)+o(T(r,G_f))=m_1T(r,f(z))+o(T(r,f(z))).
\end{align}
Applying Lemma \ref{L.5} to (\ref{pds2}), we arrive at \eqref{sm.1}. To complete the proof, we divide our analysis into the following cases:
\smallskip

\noindent
{\bf Case 1.} Let $m_1=m_2$. Then (\ref{sm.1}) gives $m_1=m_2=2$ and so from \eqref{pds2}, we obtain $G_f^2(z)+(P_1(z)f(z))^2=P(z)$, i.e.,
\begin{align}\label{sm.2.4} 
	(G_f(z)+\iota P_1(z) f(z))(G_f(z)-\iota P_1(z) f(z))=P(z).
\end{align}	
Therefore we may assume that 
\begin{align}\label{sm.2.5}
	\begin{cases}
		G_f(z)+\iota P_1(z) f(z)=\exp(\iota Q^{*}(z))\\
		G_f(z)-\iota P_1(z) f(z)=P(z)\exp(-\iota Q^{*}(z)),
	\end{cases}
\end{align}
or 
\begin{align}\label{sm.2.6}
	\begin{cases}
		G_f(z)-\iota P_1(z) f(z)=\exp(\iota Q^{*}(z))\\
		G_f(z)+\iota P_1(z) f(z)=P(z)\exp(-\iota Q^{*}(z)),\end{cases}
\end{align}
where $Q^{*}(z)$ is an entire function in $\mathbb{C}^m$. Solving, we get
\begin{align}\label{sm.2.7} 
	G_f(z)=P_2(z)f(z+c)+\sideset{}{_{j=1}^m}{\sum} Q_j(z)\partial_{z_j}(f(z))=\frac{\exp(\iota Q^{*}(z))+P(z)\exp(-\iota Q^{*}(z))}{2}
\end{align}
and 
\begin{align} \label{sm.2.8}
	f(z)=k\frac{\exp(\iota Q^{*}(z))-P(z)\exp(-\iota Q^{*}(z))}{2\iota P_1(z)},
\end{align}
where $k=1$ for \eqref{sm.2.7} and $k=-1$ for \eqref{sm.2.8}.
Using Lemma \ref{L.4e}, we find that $Q^{*}(z)$ is a non-constant polynomial in $\mathbb{C}^m$. Again from (\ref{sm.2.8}), we have	
\begin{align*} 
	\partial_{z_j}(f(z))=&k\left(\frac{1}{2P_1(z)}\partial_{z_j}( Q^{*}(z))-\frac{1}{2 \iota {P_1}^2(z)}\partial_{z_j}(P_1(z))\right)\exp (\iota Q^{*}(z))\\+&k\left(\frac{ P(z)}{2P_1(z)}\partial_{z_j}(Q^*(z))+\frac{P(z)}{2\iota {P_1}^2(z)}\partial_{z_j}(P_1(z))-\frac{1}{2 \iota P_1(z)}\partial_{z_j}( P(z))\right) \exp(-\iota Q^{*}(z)),
\end{align*}
where $j=1,2,\ldots,m$ and so 
\begin{align}\label{sm.2.9}
	G_f(z)=&P_2(z)\ol{f}(z)+\sideset{}{_{j=1}^m}{\sum} Q_j(z)\partial_{z_j}(f(z))\\=&\frac{kP_2}{2\iota \ol{P_1}}\exp(-\iota {\ol{Q^{*}}})+\nonumber k\sideset{}{_{j=1}^m}{\sum}Q_j\left(\frac{1}{2{P_1}}\partial_{z_j}({Q^*})-\frac{1}{2 \iota {P_1}^2}\partial_{z_j}({P_1})\right)\exp(\iota {Q^{*}})\\&\nonumber-\frac{kP_2{\ol P}}{2 \iota {\ol{P_1}}}\exp(-\iota {\ol{Q^{*}}})+k\sideset{}{_{j=1}^m}{\sum}Q_j\left(\frac{{P}}{2P_1}\partial_{z_j}( {Q^*})+\frac{P \partial_{z_j}({P_1})}{2\iota {P_1}^2}-\frac{\partial_{z_j}({P})}{2 \iota {P_1}}\right)\exp(-\iota {Q^{*}}).\nonumber
\end{align}
Clearly from (\ref{sm.2.7}) and (\ref{sm.2.9}), we get
\begin{align}\label{sm.2.10} 
	\alpha \exp(\iota \ol{Q^{*}}(z))+\beta \exp(-\iota \ol{Q^{*}}(z))+\gamma\exp(\iota Q^{*}(z))+\delta\exp(-\iota Q^{*}(z))\equiv 0, \end{align} 
where 
\begin{align*}
	\alpha\equiv \frac{kP_2}{2\iota \ol{P_1}}, \;\; \beta \equiv -\frac{kP_2{\ol P}}{2 \iota {\ol{P_1}}},\;\;
	\gamma \equiv k \sideset{}{_{j=1}^m}{\sum}Q_j\left(\frac{1}{2{P_1}}\partial_{z_j}({Q^*})-\frac{1}{2 \iota {P_1}^2}\partial_{z_j}({P_1})\right)-1/2
\end{align*}
and
\begin{align*}
	\delta \equiv k\sideset{}{_{j=1}^m}{\sum}Q_j\left(\frac{{P}}{2P_1}\partial_{z_j}( {Q^*})+\frac{{P}}{2\iota {P_1}^2}\partial_{z_j}({P_1})-\frac{1}{2 \iota {P_1}}\partial_{z_j}( {P})\right)-\frac{P}{2}.
\end{align*}
Since, $kPP_1P_2\not\equiv 0$, it follows that $\alpha\beta \not\equiv 0.$ Now, from (\ref{sm.2.10}), we find
\begin{align}\label{sm.2.11}
(-\beta/\alpha) \exp (-2\iota \ol{Q^*})-(\gamma/\alpha) \exp (\iota ({Q^*}-\ol{Q^*})) -(\delta/\alpha)\exp  (-\iota({Q^*}+\ol{Q^*}))\equiv 1.
\end{align}
Since $Q^*(z)$ is a non-constant polynomial in $\mathbb{C}^m$, we deduce that $\ol{Q^*}+Q^*$ is also a non-constant polynomial in $\mathbb{C}^m$. 

\noindent
We consider the following four sub-cases.
\medskip

\noindent
{\bf Sub-case 1.1.} Assume that $\gamma\equiv \delta \equiv 0.$ Then (\ref{sm.2.11}) reduces to $-(\beta/\alpha)\exp  (-2\iota \ol{Q^*})\equiv 1$, which is impossible by Lemma \ref{L.3}.
\medskip

\noindent
{\bf Sub-case 1.2.} Assume that $\gamma\equiv0$ and $\delta \not\equiv 0.$ Then (\ref{sm.2.11}) reduces to 
\begin{align}\label{sm.2.12}
	-(\beta/\alpha) \exp (-2\iota \ol{Q^*})-(\delta/\alpha)\exp  (-\iota({Q^*}+\ol{Q^*}))\equiv 1.
\end{align} 
It is easy to verify that $T(r,\beta/\alpha)=o(T(r, \exp  (-2\iota Q^*)))$, $T(r,\beta/\alpha)=o(T(r, \exp (-\iota (\ol{Q^*}+Q^*))))$, $T(r,\delta/\alpha)=o(T(r,\exp Q^*))$ and $T(r,\delta/\alpha)=o(T(r, \exp (-\iota (\ol{Q^*}+Q^*)))).$ 

\noindent
Let $F=-\frac{\beta}{\alpha} \exp (-2\iota \ol{Q^*})$ and $G=-\frac{\delta}{\alpha}\exp  (-\iota({Q^*}+\ol{Q^*}))$. Using Lemma \ref{L.2}, we get
\begin{align*}
	T(r,F) \le& \ol N(r, F)+\ol N(r,0;F)+\ol N(r,1;F) +o(T(r,F))\\=& \ol N(r, F)+\ol N(r,0;F)+\ol N(r,0;G) +o(T(r,F))\\\le & \ol N(r,\beta/\alpha)+\ol N(r,\alpha/\beta)+\ol N(r,\alpha/\delta) +o(T(r,F)) \le  o(T(r,F)),
\end{align*}
which leads to a contradiction.
\medskip

\noindent
{\bf Sub-case 1.3.} Assume that $\gamma\not\equiv0$ and $\delta \equiv 0$. Then (\ref{sm.2.11}) reduces to
\begin{align}\label{sm.2.13}
	(-\beta/\alpha) \exp (-2\iota \ol{Q^*})-(\gamma/\alpha) \exp (\iota ({Q^*}-\ol{Q^*})) \equiv 1.
\end{align}
If ${Q^*}-\ol{Q^*}$ is constant, then using Lemma \ref{L.3} to \eqref{sm.2.13}, we get a contradiction. Hence ${Q^*}-\ol{Q^*}$ is non-constant.  Now, proceeding in the same way as done in the proof of Sub-case 1.2, we get a contradiction.
\medskip

\noindent
{\bf Sub-case 1.4.} Assume that $\gamma\not\equiv0$ and $\delta\not \equiv 0$. Then applying Lemma \ref{L.6} to (\ref{sm.2.11}), we get 
\begin{align}\label{sm.2.14}
	-(\gamma/\alpha) \exp (\iota ({Q^*}-\ol{Q^*})) \equiv 1,
\end{align}
which shows that ${Q^*}-\ol{Q^*}$ is a constant and so we may assume that
\begin{align}\label{sm.2.15}
	Q^*(z)=d_1z_1+d_2z_2+\ldots+d_mz_m+\tilde Q^*(z),
\end{align} 
where $d_1,d_2,\ldots,d_m\in\mathbb{C}$ and $\tilde Q^*(z)$ is a polynomial in $\mathbb{C}^m$ such that $\tilde Q^*(z+c)\equiv \tilde Q^*(z)$.

\noindent	
Substituting \eqref{sm.2.14} into \eqref{sm.2.11}, we get
\begin{align}\label{sm.2.16}
	\exp (-\iota (\ol{Q^*}-Q^*))\equiv -\beta/\delta.
\end{align}	
Let $d^*=\exp (-\iota (\ol{Q^*}-Q^*))$. Now from \eqref{sm.2.14} and \eqref{sm.2.16}, we get $d^*\equiv -\frac{\beta}{\delta} \equiv - \frac{\gamma}{\alpha}$ and so $\beta\alpha \equiv\gamma \delta$. Therefore, we obtain
\begin{align}\label{sm.2.17}
	P_2d^*+\ol {P_1} \sideset{}{_{j=1}^m}{\sum} Q_j\left(\frac{\iota}{P_1} \partial_{z_j}( Q^*)-\frac{1}{P_1^2}\partial_{z_j}( P_1)\right)\equiv \frac{\iota  \ol{P_1}}{k}.
\end{align} 

Finally, from \eqref{sm.2.8} we have 
\begin{align*}
	f(z)=k\frac{\exp(\iota Q^{*}(z))-P(z)\exp(-\iota Q^{*}(z))}{2\iota P_1(z)},	
\end{align*} where $Q^*(z)=d_1z_1+d_2z_2+\ldots+d_mz_m+\tilde Q^*(z),$
where $d_1,d_2,\ldots,d_m\in\mathbb{C}$ and $\tilde Q^*(z)$ is a polynomial in $\mathbb{C}^m$ such that $\tilde Q^*(z+c)\equiv\tilde  Q^*(z)$, \begin{align*}
	P_2e^{-\iota (d_1c_1+d_2c_2+\ldots+d_mc_m)}+\ol {P_1} \sideset{}{_{j=1}^m}{\sum} Q_j\left(\frac{\iota}{P_1} \partial_{z_j}( Q^*)-\frac{1}{P_1^2}\partial_{z_j}( P_1)\right)\equiv\frac{\iota  \ol{P_1}}{k}.
\end{align*} 
\smallskip

\noindent
{\bf Case 2.} Let $m_1\neq m_2$. We now consider the following two sub-cases.
\smallskip

\noindent
{\bf Sub-case 2.1.} Suppose $m_1 >m_2$. Then \eqref{sm.1} yields $m_2=1$ and so from \eqref{pds2}, we have 
\begin{align}\label{sm.2.19}
	(P_1f(z))^{m_1}+G_f(z)=P(z).
\end{align}
Now using Lemmas \ref{L.1}, \ref{L.3}, \ref{L.8} and (\ref{sm.11}) to (\ref{sm.2.19}), we get
\begin{align*}
	\parallel\;m_1 T(r, f))=&T\left(r, \left(P_1(z)f(z)\right)^{m_1}\right)+o(T(r, f))\nonumber\\
	=&T\left(r,G_f(z)-P(z)\right)+o(T(r, f)) \nonumber\\
	\leq&m\left(r,G_f(z)\right)+o(T(r, f)) \nonumber\\
	\leq &m\left(r, \sideset{}{_{j=1}^m}{\sum} Q_j\partial_{z_j}(f(z))/f(z)\right)+m\left(r,f(z)/\ol{f}(z)\right)+m(r, \ol{f}(z))+o(T(r, f))\nonumber \\
	=&T(r, \ol{f}(z))+o(T(r, f))=T(r,f)+o(T(r, f))
\end{align*}
and so $\parallel (m_1-1) T(r, f) \leq o(T(r, f))$. Since $m_1\geq 2$, we get a contradiction.
\smallskip

\noindent
{\bf Sub-case 2.2.} Suppose $m_2>m_1$, $\sum_{j=1}^m \partial_{z_j}(\ol{P_1}(z))=P_2(z)$, $Q_i(z)=P_1(z)$ for $j=1,2,\ldots,m$ and $\sum_{j=1}^m Q_j(z)\partial_{z_j}( f(z))$ is $c$-periodic. From \eqref{sm.1}, we have $m_1=1$ and so (\ref{pds2}) gives
\begin{align}\label{sm.2.21} 
	P_1(z)f(z)+G_f^{m_2}(z)=P(z),
\end{align}
Differentiating (\ref{sm.2.21}) partially with respect to $z_j$, we get
\begin{align*} 
m_2G_f^{m_2-1}(z)\partial_{z_j}(G_f(z))+\partial_{z_j}(P_1(z))f(z)+P_1(z)\partial_{z_j}(f(z))=\partial_{z_j}( P(z)),
\end{align*}
i.e.,
\begin{align}\label{sm.2.22} 
&m_2G_f^{m_2-1}\sideset{}{_{j=1}^m}{\sum} \partial_{z_j}(G_f(z))+
\sideset{}{_{j=1}^m}{\sum} \partial_{z_j}(P_1(z)) f(z)+\sideset{}{_{j=1}^m}{\sum} P_1\partial_{z_j}(f(z))
\\=&\sideset{}{_{j=1}^m}{\sum}\partial_{z_j}(P(z)).\nonumber
\end{align}
Now in view of \eqref{sm.2.1}, we get from \eqref{sm.2.22} that
\begin{align}\label{sm.2.23}
m_2\ol{G_f}^{m_2-1}(z)\tilde \beta(z)+G_f(z)=\tilde \beta_1(z),
\end{align}
where
\begin{align}\label{sm.2.24}
\tilde\beta(z)=\sideset{}{_{j=1}^m}{\sum}\partial_{z_j}(\ol{G_f}(z))
\end{align}
and
\begin{align}\label{sm.2.25}
\tilde\beta_1(z)=\sideset{}{_{j=1}^m}{\sum}\partial_{z_j}(\ol{P}(z)).
\end{align}
If $\tilde\beta\equiv 0$, then from (\ref{sm.2.23}), we see that 
$G_f(z)$ is a polynomial in $\mathbb{C}^m$. Therefore from \eqref{sm.2.21}, we find that $f(z)$ is a polynomial in $\mathbb{C}^m$, which is impossible. 
Hence $\tilde\beta\not\equiv 0$.
We know that $m_2\geq 2$. Using Lemma \ref{L.4} to (\ref{sm.2.23}), we get $\parallel m (r,\tilde\beta(z))=o(T(r, \ol{G_f}(z))).$
Since $\tilde\beta$ is an entire function, we have
\begin{align}\label{sm.2.26}
\parallel T(r,\tilde\beta(z))=o(T(r, \ol{G_f}(z))).
\end{align}

\noindent
Based on the above estimations, we now consider the following sub-cases:
\smallskip

\noindent
{\bf Sub-case 2.2.1.} Let $m_2-1\geq 2$. Now using Lemma \ref{L.4} to (\ref{sm.2.23}), we get $\parallel m(r,\tilde \beta \ol{G_f}(z))=o(T(r,\ol{G_f}(z)))$.
Since $\ol{G_f}(z)$ is an entire function and $\parallel T(r,\tilde \beta(z))=o(T(r,\ol{G_f}(z)))$, we have 
\begin{align*}
	\parallel\;T(r,\ol{G_f}(z))=m(r,\ol{G_f}(z))\leq m(r,\tilde\beta \ol{G_f}(z))+m(r,1/\tilde\beta(z))\leq o(T(r,\ol{G_f}(z)))
\end{align*}
and so in view of Lemma \ref{L.7} and using (\ref{sm.2.2}), we have $\parallel T(r,f)=o(T(r,f))$, which is impossible.\par
\smallskip

\noindent
{\bf Sub-case 2.2.2.} Let $m_2-1=1$. Then from (\ref{sm.2.23}), we get
\begin{align}\label{sm.2.27} 
	2\ol{G_f}(z)\tilde\beta(z)+G_f(z)=\tilde\beta_1(z).
\end{align}
Differentiating (\ref{sm.2.27}) partially with respect to $z_j$, we get
\begin{align*}
	2\partial_{z_j}(\ol{G_f}(z))\tilde\beta(z)+2\ol{G_f}(z)\partial_{z_j}(\tilde\beta(z))+\partial_{z_j}( G_f(z))=0,
\end{align*}
i.e.,
\begin{align}\label{sm.2.28}
	2\tilde\beta(z)\sideset{}{_{j=1}^m}{\sum} \partial_{z_j}(\ol{G_f}(z))+2\ol{G_f}(z)\sideset{}{_{j=1}^m}{\sum}\partial_{z_j}( \tilde\beta(z))+\sideset{}{_{j=1}^m}{\sum} \partial_{z_j}(G_f(z))=0.
\end{align}
Using \eqref{sm.2.24} to \eqref{sm.2.28}, we get
\begin{align}\label{sm.2.29} 
	2\tilde\beta^2(z)+2\ol{G_f}(z)\sideset{}{_{j=1}^m}{\sum}\partial_{z_j}( \tilde\beta(z))+\ul{\tilde \beta}(z)=0.
\end{align}
\medskip

\noindent
Suppose $\sum_{j=1}^m\partial_{z_j}(\tilde\beta(z))\not\equiv 0$. Now in view of (\ref{sm.2.26}) and using Lemma \ref{L.3} to (\ref{sm.2.29}), we get  $\parallel T(r, \ol{G_f}(z))=o(T(r, \ol{G_f}(z)))$ and so from Lemma \ref{L.7} and (\ref{sm.2.2}), we get 
$\parallel T(r,f)=o(T(r,f)$, which is impossible. Hence $\sum_{j=1}^m\partial_{z_j}(\tilde\beta(z))\equiv 0$. Then from (\ref{sm.2.29}), we have
\begin{align}\label{sm.2.30} 
	2\tilde \beta^2(z)=-\ul{\tilde \beta}(z).
\end{align}
It is easy to prove from (\ref{sm.2.30}) that $\tilde \beta(z)$  is a constant, say $k_1(\neq 0)$. Then from (\ref{sm.2.30}), we have 
$2k_1^2=-k$, which implies that $k_1=-1/2$. Therefore 
\begin{align}\label{sm.2.31} 
	\tilde\beta(z)=\sideset{}{_{j=1}^m}{\sum}\partial_{z_j}( \ol{G_f}(z))=-1/2.
\end{align} 
Again from (\ref{sm.2.25}) and \eqref{sm.2.27}, we get 
\begin{align}\label{sm.2.32}
	\ol{G_f}(z)=G_f(z)-\sideset{}{_{j=1}^m}{\sum}\partial_{z_j}(\ol{P}(z)).
\end{align}
The Lagrange's auxiliary equations corresponding to (\ref{sm.2.31}) are as follows
\begin{align}\label{sm.2.33} 
	\displaystyle
	\frac{d z_1}{1}=\frac{dz_2}{1}=\ldots=\frac{dz_m}{1}=\frac{dG_f(z)}{-1/2}=\frac{dz_1+dz_2+\ldots+dz_m}{m}.
\end{align}
Clearly $z_2-z_1=d_1$, $z_3-z_2=d_2$, $\ldots$, $z_m-z_{m-1}=d_{m-1}$ and $mG_f(z)+\frac{1}{2}(z_1+z_2+\ldots+z_m)=d_m$ be independent solutions of (\ref{sm.2.31}). Therefore the general solution of (\ref{sm.2.31}) is as follows
\begin{align}\label{sm.2.34} 
	\displaystyle
	G_f(z)=-(z_1+z_2+\ldots+z_m)/2m+g_1(z),
\end{align}
where $g_1(z)=g_1(z_2-z_1,z_3-z_2,\ldots,z_m-z_{m-1})$ is an arbitrary finite-order transcendental entire function. Using \eqref{sm.2.34} to (\ref{sm.2.32}), we have
\begin{align*}
	\displaystyle
	&-(z_1 + c_1+z_2+c_2+\ldots+z_m+c_m)/2m + g_1(z+c)+\sideset{}{_{j=1}^m}{\sum}\partial_{z_j}(\ol{P}(z))\\=&-(z_1+z_2+\ldots+z_m)/2m + g_1(z),
\end{align*}
which implies
\begin{align}\label{sm.2.35} 
	\displaystyle
	g_1(z+c)-g_1(z)=(c_1+c_2+\ldots+c_m)/2m-\sideset{}{_{j=1}^m}{\sum}\partial_{z_j}(\ol{P}(z)).
\end{align}
Finally from \eqref{sm.12}, we have 
\begin{align*}
	f(z)=\frac{1}{P_1(z)}\left(-\Big(-(z_1+z_2+\ldots+z_m)/2m+g_1(z)\Big)^2+P(z)\right),
\end{align*}
where $g_1(z)=g_1(z_2-z_1,z_3-z_2,\ldots,z_m-z_{m-1})$ is an arbitrary finite-order transcendental entire function such that (\ref{sm.2.35}) holds.
\end{proof}

\section {{\bf The transcendental entire solutions of Eq. \eqref{pds3}}} \label{sec5}

\noindent
This section is devoted to stating our main result, which serves as a higher-dimensional version of Theorem C.
\begin{theo}\label{t2.3} Let $P(\not\equiv 0)$, $P_1$, $P_2$ and $Q_i\;(i=1,2,\ldots,m)$ are polynomials in $\mathbb{C}^m$ such that $P_2Q_j\not\equiv 0\;(j=1,2,\ldots,m)$ and $P(z)$ is irreducible in $\mathbb{C}^m$. Suppose $f(z)$ is a finite-order transcendental entire solution to Eq. (\ref{pds3}) such that $\sum_{j=1}^m Q_j(z)\partial_{z_j}( f(z))$ is transcendental. Then
\begin{enumerate}
\item[\em{(1)}] if $m_1=m_2$, $P_1=MP_2$, where $M$ is a constant and $Q_j=Q$ for $j=1,2,\ldots,m$, then $m_1=m_2=2$. Also $f(z)$ satisfies
\begin{align*}
\sideset{}{_{j=1}^{m}}{\sum}\partial_{z_j}( f(z))=\frac{\exp(\iota Q^{*}(z))+P(z)\exp(-\iota Q^{*}(z))}{2Q(z)}
\end{align*}
and
 \begin{align*}
	Mf(z+c)+f(z)\equiv k\frac{\exp(\iota Q^{*}(z))-P(z)\exp(-\iota Q^{*}(z))}{2\iota P_2(z)}\;(k=\pm 1),
	\end{align*}
where $Q^*$ satisfies the following two cases:
\begin{enumerate}
\item[\em{(1A)}] whenever $M=0$, then $Q^*$ is a non-constant polynomial in $\mathbb{C}^m$ such that 
\begin{align*}
	kQ\sideset{}{_{j=1}^m}{\sum}\left(\frac{1}{P_2}\partial_{z_j}( Q^*)+\frac{\iota}{P_2^2}\partial_{z_j}( P_2)\right)\equiv 1\;\;
\text{and}\;\;
	\frac{1}{P}\sideset{}{_{j=1}^m}{\sum}\partial_{z_j}(P)\equiv \frac{2}{P_2}\sideset{}{_{j=1}^m}{\sum}\partial_{z_j}(P_2);
\end{align*}
\item[\em{(1B)}] whenever $M\ne 0$, then $Q^*(z)=d_1z_1+d_2z_2+\ldots+d_mz_m+\tilde Q^*(z),$ $d_1,d_2,\ldots,d_m\in\mathbb{C}$, $\tilde Q^*(z)$ is a polynomial in $\mathbb{C}^m$ such that $\tilde Q^*(z+c)\equiv \tilde Q^*(z)$ and
\begin{align*}
	Me^{-\iota (d_1c_1+d_2c_2+\ldots+d_mc_m)}+\frac{\ol Q}{Q}\equiv k\frac{\ol Q}{P_2}\sideset{}{_{j=1}^{m}}{\sum}\partial_{z_j}(Q^*)+k\frac{\iota \ol Q}{P_2^2}\sideset{}{_{j=1}^{m}}{\sum}\partial_{z_j}(P_2);
\end{align*}
\end{enumerate}
\smallskip 
		
\item[\em{(ii)}] if $m_1>m_2$, $P_1(z)=P_2(z)\equiv 1$, $Q_j(z)\equiv 1\;(j=1,2,\ldots,m)$ and $\sum_{j=1}^m\partial_{z_j}(f(z))$ is $c$-periodic, then $m_1=2$, $m_2=1$ and $f(z)$ satisfies 
\begin{align*}
f(z+c)+f(z)=-\left(-(z_1+z_2+\ldots+z_m)/2m+g_1(z)\right)^2+P(z)
\end{align*}
and
\begin{align*}
	\sideset{}{_{j=1}^m}{\sum}\partial_{z_j}(f(z))=-(z_1+z_2+\ldots+z_m)/m+g_1(z),
\end{align*}
where $g_1(z)=g_1(z_2-z_1,z_3-z_2,\ldots,z_m-z_{m-1})$ is an arbitrary finite-order transcendental entire function and $\sum_{j=1}^m\partial_{z_j}(P(z))=0$.		
\item[(iii)] if $m_2>m_1$, $P_1(z)=P_2(z)\equiv 1$, $Q_j(z)\equiv 1\;(j=1,2,\ldots,m)$ and $\sum_{j=1}^m\partial_{z_j}(f(z))$ is $c$-periodic, then $f(z)$ can not be a solution to Eq. (\ref{pds3}).
\end{enumerate}
\end{theo}
In the special case where $P_1=0$, $P_2=P=1$, $Q=1$ and $m_1=m_2=2$, Theorem \ref{t2.3} reduces to the following corollary, refining \cite[Corollary 1.4]{XC2}.
\begin{cor} All finite-order transcendental entire solutions of 
\begin{align*}
	\left(\sideset{}{_{j=1}^m}{\sum}\partial_{z_j}(f(z))\right)^2+f^2(z)=1
\end{align*}
must satisfy	
\begin{align*}
f(z)=\frac{\exp(\iota Q^{*}(z))-\exp(-\iota Q^{*}(z))}{2\iota},
\end{align*}
where $Q^*$ is a non-constant polynomial in $\mathbb{C}^m$ such that $\sideset{}{_{j=1}^m}{\sum}\partial_{z_j}(Q^*)\equiv 1$.	
\end{cor}
Setting $P_1=1$ $P_2=-1$, $P=1$, $Q=1$ and $m_1=m_2=2$ in Theorem \ref{t2.3} yields the following corollary, which improves upon \cite[Theorem 1.4]{Haldar_Mediterian_2023}.
\begin{cor} All finite-order transcendental entire solutions of 
\begin{align*}
\left(\sideset{}{_{j=1}^m}{\sum}\partial_{z_j}(f(z))\right)^2+\left(f(z+c)-f(z)\right)^2=1
\end{align*}
must satisfy	
\begin{align*}
\sideset{}{_{j=1}^{m}}{\sum}\partial_{z_j}(f(z))=\frac{\exp(\iota Q^{*})+\exp(-\iota Q^{*})}{2}\;
\text{and}\;
f(z+c)-f(z)=\frac{\exp(\iota Q^{*})-\exp(-\iota Q^{*})}{2\iota},
\end{align*}
where $Q^*(z)=d_1z_1+d_2z_2+\ldots+d_mz_m+\tilde Q^*(z),$, $d_1,d_2,\ldots,d_m\in\mathbb{C}$, $\tilde Q^*(z)$ is a polynomial in $\mathbb{C}^m$ such that $\tilde Q^*(z+c)\equiv \tilde Q^*(z)$ and
\begin{align*}
\exp({-\iota (d_1c_1+d_2c_2+\ldots+d_mc_m)})=\sideset{}{_{j=1}^{m}}{\sum}\partial_{z_j}(Q^*)+1.
\end{align*}	
\end{cor}

\begin{proof}[{\bf Proof of Theorem \ref{t2.3}}] 
Let $f(z)$ be any finite-order transcendental entire solution to the Eq. (\ref{pds3}) such that $\sum_{j=1}^m Q_j(z)\partial_{z_j}(f(z))$ is also transcendental. Take
\begin{align}\label{sm.5.1}
	H_f(z)=\sideset{}{_{j=1}^m}{\sum} Q_j(z)\partial_{z_j}(f(z)) \;\;\text{and}\;\; I_f(z)=P_1(z)f(z+c)+P_2(z)f(z).
\end{align}
Therefore \eqref{pds3} gives 
\begin{align}\label{sm.5.2}
(H_f(z))^{m_1}+(I_f(z))^{m_2}=P(z).
\end{align}
Using Lemma \ref{L.3} to (\ref{pds3}), we get
\begin{align}\label{sm.5.3} \parallel\;m_1T(r,H_f)+o(T(r,H_f))=m_2T(r,I_f)+o(T(r,f(z))).
\end{align}
Applying Lemma \ref{L.5} to (\ref{pds3}), we arrive at \eqref{sm.1}. We now consider the following cases.
\smallskip

\noindent
{\bf Case 1.} Let $m_1=m_2$, $P_1=MP_2$, where $M$ is a constant and $Q_j=Q$ for $j=1,2,\ldots,m$. Then (\ref{sm.1}) gives $m_1=m_2=2$ and so from \eqref{sm.5.2}, we get $H_f^2(z)+I_f^2(z)=P(z)$, i.e.,
\begin{align}\label{sm.5.5} 
(H_f(z)+\iota I_f(z))(H_f(z)-\iota I_f(z))=P(z).
\end{align}
Therefore we may assume that 
\begin{align}\label{sm.5.6}
\begin{cases}
H_f(z)+\iota I_f(z)=\exp(\iota Q^{*}(z))\\
H_f(z)-\iota I_f(z)=P(z)\exp(-\iota Q^{*}(z)),
\end{cases}
\end{align}
or 
\begin{align}\label{sm.5.7}
\begin{cases}
H_f(z)-\iota I_f(z)=\exp(\iota Q^{*}(z))\\
H_f(z)+\iota I_f(z)=P(z)\exp(-\iota Q^{*}(z)),\end{cases}
\end{align}
where $Q^{*}(z)$ is an entire function in $\mathbb{C}^m$. Since, $Q_1(z)=Q_2(z)=Q_3(z)=\ldots=Q_m(z)=Q(z)$ and $P_1(z)=MP_2(z)$, solving \eqref{sm.5.6} and \eqref{sm.5.7} we get
\begin{align}\label{sm.5.8} 
	H_f(z)=Q(z)\sideset{}{_{j=1}^m}{\sum}\partial_{z_j}(f(z))=\frac{\exp(\iota Q^{*}(z))+P(z)\exp(-\iota Q^{*}(z))}{2}
\end{align}
and 
\begin{align} \label{sm.5.9}
Mf(z+c)+f(z)=k\frac{\exp(\iota Q^{*}(z))-P(z)\exp(-\iota Q^{*}(z))}{2\iota P_2(z)},
\end{align}
where $k=1$ for \eqref{sm.5.5} and $k=-1$ for \eqref{sm.5.6}.
Now by Lemma \eqref{L.4e}, we deduce that $Q^{*}(z)$ is a non-constant polynomial in $\mathbb{C}^m$.  Using \eqref{sm.5.8} and \eqref{sm.5.9} and performing a simple computation, we obtain
\begin{align}\label{sm.5.10}
(M/\ol Q) \exp{(\iota \ol{Q}^{*})}+(M\ol P/\ol Q) \exp{(-\iota \ol{Q}^{*})}+\alpha \exp{(\iota {Q}^{*})}+\beta \exp{(-\iota {Q}^{*})} \equiv 0,
\end{align}
where 
\begin{align*}\alpha \equiv 1/Q- k\sideset{}{_{j=1}^m}{\sum}\left(\frac{1}{P_2}\partial_{z_j}(Q^*)+\frac{\iota}{P_2^2}\partial_{z_j}(P_2)\right)
\end{align*}
and
\begin{align*}
	\beta \equiv P/Q-k \sideset{}{_{j=1}^m}{\sum} \left(\frac{P}{P_2}\partial_{z_j}(Q^*)+\frac{\iota}{P_2}\partial_{z_j}(P)-(\iota P/P_2^2)\partial_{z_j}( P_2)\right).
\end{align*}
We now consider the following sub-cases:
\smallskip

\noindent
{\bf Sub-case 1.1.} Let $M=0.$ Then from \eqref{sm.5.10}, we have $\alpha \exp {(2\iota Q^*)}\equiv -\beta.$ Since $Q^*$ is a non-constant polynomial, we have $\alpha \equiv \beta \equiv 0.$ Therefore, we have 
\begin{align*}
kQ\sideset{}{_{j=1}^m}{\sum}\left(\frac{1}{P_2}\partial_{z_j}(Q^*)+\frac{\iota}{P_2^2}\partial_{z_j}(P_2)\right)\equiv 1\;\;
\text{and}\;\;
\frac{1}{P}\sideset{}{_{j=1}^m}{\sum}\partial_{z_j}(P)\equiv \frac{2}{P_2}\sideset{}{_{j=1}^m}{\sum}\partial_{z_j}(P_2).
\end{align*}
\smallskip

\noindent
{\bf Sub-case 1.2.} Let $M \ne 0.$ From (\ref{sm.5.10}), we find
\begin{align}\label{sm.5.11}
-\ol P \exp (-2\iota \ol{Q^*})-(\alpha\ol Q/M) \exp (\iota ({Q^*}-\ol{Q^*})) -(\beta\ol Q/M)\exp  (-\iota({Q^*}+\ol{Q^*}))\equiv 1.
\end{align}
Since $Q^*(z)$ is a non-constant polynomial in $\mathbb{C}^m$, we deduce that $\ol{Q^*}+Q^*$ is also a non-constant polynomial in $\mathbb{C}^m$. 
We consider the following four sub-cases.
\medskip

\noindent
{\bf Sub-case 1.2.1.} Assume that $\alpha\equiv \beta \equiv 0.$ Then (\ref{sm.5.11}) reduces to $-\ol P\exp  (-2\iota \ol{Q^*})\equiv 1$, which is impossible by Lemma \ref{L.3}.
\medskip

\noindent
{\bf Sub-case 1.2.2.} Assume that $\alpha\equiv0$ and $\beta \not\equiv 0.$ Then (\ref{sm.5.11}) reduces to 
\begin{align*}
	-\ol P \exp (-2\iota \ol{Q^*}) -(\beta\ol Q/M)\exp  (-\iota({Q^*}+\ol{Q^*}))\equiv 1.
\end{align*} 
Proceeding analogously to the proof of \textit{Sub-Case 1.2} of \textit{Theorem \ref{t2.2}}, we arrive at a contradiction
\medskip

\noindent
{\bf Sub-case 1.2.3.} Assume that $\alpha\not\equiv0$ and $\beta \equiv 0$. Then (\ref{sm.5.11}) reduces to
\begin{align*}
	-\ol P \exp (-2\iota \ol{Q^*})-(\alpha\ol Q/M) \exp (\iota ({Q^*}-\ol{Q^*})) \equiv 1.
\end{align*}
Following the argument of \textit{Sub-Case 1.3} in \textit{Theorem \ref{t2.2}}, we obtain a contradiction.
\medskip

\noindent
{\bf Sub-case 1.2.4.} Assume that $\alpha\not\equiv0$ and $\beta\not \equiv 0$. Applying Lemma \ref{L.6} to (\ref{sm.5.11}), we get 
\begin{align}\label{sm.5.12}
	-(\alpha\ol Q/M) \exp (\iota ({Q^*}-\ol{Q^*})) \equiv 1,
\end{align}
which shows that ${Q^*}-\ol{Q^*}$ is a constant and so similar to {\it Theorem \ref{t2.2}}, we may assume that
\begin{align}\label{sm.5.13}
	Q^*(z)=d_1z_1+d_2z_2+\ldots+d_mz_m+\tilde Q^*(z),
\end{align} 
where $d_1,d_2,\ldots,d_m\in\mathbb{C}$ and $\tilde Q^*(z)$ is a polynomial in $\mathbb{C}^m$ such that $\tilde Q^*(z+c)\equiv \tilde Q^*(z)$.	
Substituting \eqref{sm.5.12} into \eqref{sm.5.11}, we get
\begin{align}\label{sm.5.14}
	\exp (-\iota (\ol{Q^*}-Q^*))\equiv -\beta\ol Q/M\ol P.
\end{align}	

Let $d^*=\exp (-\iota (\ol{Q^*}-Q^*))$. Now from \eqref{sm.5.12} and \eqref{sm.5.14}, we get $d^*\equiv -\beta\ol Q/M\ol P \equiv - \alpha\ol Q/M$ and so $\beta \equiv\ol P\alpha$. Therefore, we obtain
\begin{align*}
	Me^{-\iota (d_1c_1+d_2c_2+\ldots+d_mc_m)}+\frac{\ol Q}{Q}\equiv k\frac{\ol Q}{P_2}\sideset{}{_{j=1}^{m}}{\sum}\partial_{z_j}(Q^*)+k\frac{\iota \ol Q}{P_2^2}\sideset{}{_{j=1}^{m}}{\sum}\partial_{z_j}(P_2).
\end{align*}
\smallskip

\noindent
{\bf Case 2.} Let $m_1\neq m_2$. In this case, we have $P_1(z)=P_2(z)\equiv 1$, $Q_j(z)\equiv 1$ for $j=1,2,\ldots,m$ and $\sum_{j=1}^m \partial_{z_j}(f(z))$ is $c$-periodic. We now consider the following two sub-cases.
\smallskip

\noindent
{\bf Sub-case 2.1.} Suppose $m_2 >m_1$. Then \eqref{sm.1} gives $m_1=1$ and so from \eqref{pds3}, we get
\begin{align}\label{sm.321.1}
H_f(z)+(I_f(z))^{m_2}=P(z).
\end{align}
By the given conditions, we have
\begin{align}\label{sm.321.2}
2H_f(z)=2\sideset{}{_{j=1}^m}{\sum} Q_j(z)\partial_{z_j}(f(z))=2\sideset{}{_{j=1}^m}{\sum}\partial_{z_j}(f(z))=\sideset{}{_{j=1}^m}{\sum}\partial_{z_j}(I_f(z))
\end{align}
and so \eqref{sm.321.1} gives
\begin{align}\label{sm.321.3}
\frac{1}{2}\sideset{}{_{j=1}^m}{\sum}\partial_{z_j}(I_f(z))+(I_f(z))^{m_2}=P(z).
\end{align}	
Now using Lemmas \ref{L.1} and \ref{L.3} to (\ref{sm.321.3}), we get
\begin{align*}
\parallel\;m_2 T(r, I_f(z)))=&T\left(r, \left(I_f(z)\right)^{m_2}\right)+o(T(r, I_f))\nonumber\\
	=&T\left(r,\frac{1}{2}\sideset{}{_{j=1}^m}{\sum}\partial_{z_j}(I_f(z))-P(z)\right)+o(T(r, I_f)) \nonumber\\
	=&T\left(r,\frac{1}{2}\sideset{}{_{j=1}^m}{\sum}\partial_{z_j}(I_f(z))\right)+o(T(r, I_f)) \nonumber\\
	=&m\left(r,\frac{1}{2}\sideset{}{_{j=1}^m}{\sum}\partial_{z_j}(I_f(z))\right)+o(T(r, I_f)) \nonumber\\
	\leq &m\left(r, \sideset{}{_{j=1}^m}{\sum}\partial_{z_j}(I_f(z))/I_f(z)\right)+m(r, I_f(z))+o(T(r, I_f))\nonumber \\
	=&T(r, I_f)+o(T(r, I_f))
\end{align*}
and so $\parallel (m_2-1) T(r, I_f) \leq o(T(r, I_f))$. Since $m_2\geq 2$, we get a contradiction.
\smallskip

\noindent
{\bf Sub-case 2.2.} Suppose $m_1 >m_2$. From \eqref{sm.1}, we have $m_2=1$. Therefore using \eqref{sm.321.2} to \eqref{pds3}, we obtain 
\begin{align}\label{sm.322.1}
\left(\frac{1}{2}\sideset{}{_{j=1}^m}{\sum}\partial_{z_j}(I_f(z))\right)^{m_1}+I_f(z)=P(z).
\end{align}	
Differentiating (\ref{sm.322.1}) partially with respect to $z_k$, we get
\begin{align*} 
	m_1\left(\frac{1}{2}\sideset{}{_{j=1}^m}{\sum}\partial_{z_j}(I_f(z))\right)^{m_1-1}\frac{1}{2}\sideset{}{_{j=1}^m}{\sum}\partial^2_{z_kz_j}(I_f(z))+\partial_{z_k}(I_f(z))=\partial_{z_k}(P(z)),
\end{align*}
i.e.,
\begin{align}\label{sm.322.2} 
m_1\tilde \alpha^{m_1-1}(z)\tilde \beta(z)+2\tilde\alpha(z)=\tilde\beta_1(z),
\end{align}
where
\begin{align}\label{sm.322.3}
\tilde \alpha(z)=\frac{1}{2}\sideset{}{_{j=1}^m}{\sum}\partial_{z_j}(I_f(z)),	
\end{align}
\begin{align}\label{sm.322.4}
	\tilde \beta(z)=\frac{1}{2}\sideset{}{_{k=1}^m}{\sum}\sideset{}{_{j=1}^m}{\sum}\partial^2_{z_kz_j}(I_f(z))=\sideset{}{_{k=1}^m}{\sum}\partial_{z_k}(\tilde \alpha(z))
\end{align}
and
\begin{align}\label{sm.322.5}
	\tilde \beta_1(z)=\sideset{}{_{k=1}^m}{\sum}\partial_{z_k}(P(z)).
\end{align}
If possible, suppose $\tilde\beta\equiv 0$. Then from (\ref{sm.322.2}) and \eqref{sm.322.3}, we have 
\begin{align*}
2\tilde\alpha(z)=\sideset{}{_{j=1}^m}{\sum}\partial_{z_j}(I_f(z))=\sideset{}{_{k=1}^m}{\sum}\partial_{z_k}(P(z)),
\end{align*}
which shows that $H_f(z)$ is a polynomial in $\mathbb{C}^m$. Therefore, we get a contradiction. 
Hence $\tilde\beta\not\equiv 0$.
We know that $m_2\geq 2$. Using Lemma \ref{L.4} to (\ref{sm.322.2}), we get $\parallel m(r,\tilde\beta)=o(T(r, \tilde\alpha)).$
Since $\tilde\beta$ is an entire function, we have
\begin{align}\label{sm.322.6}
\parallel T(r,\tilde\beta)=o(T(r, \tilde\alpha)).
\end{align}

\noindent
We now consider the following sub-cases.
\smallskip

\noindent
{\bf Sub-case 2.2.1.} Let $m_1-1\geq 2$. Now using Lemma \ref{L.4} to (\ref{sm.322.2}), we get $\parallel m(r,\tilde\alpha\beta)=o(T(r,\tilde\alpha))$.
Since $\tilde\alpha$ is an entire function and $\parallel T(r,\tilde\beta)=o(T(r,\tilde\alpha))$, we have 
\begin{align*}
\parallel\;T(r,\tilde\alpha)=m(r,\alpha)\leq m(r,\tilde\alpha\tilde\beta)+m(r,1/\tilde\beta)\leq o(T(r,\tilde\alpha)),
\end{align*}
which contradicts the fact that $\tilde\alpha$ is non-constant.
\smallskip

\noindent
{\bf Sub-case 2.2.2.} Let $m_1-1=1$. Then from (\ref{sm.322.2}), we get
\begin{align}\label{sm.322.7} 
2\tilde\alpha(z)\tilde\beta(z)+2\tilde\alpha(z)=\tilde\beta_1(z).
\end{align}
Differentiating (\ref{sm.322.7}) partially with respect to $z_k$, we get
\begin{align*}
\partial_{z_k}(\tilde\alpha(z))\tilde\beta(z)+\tilde\alpha(z)\partial_{z_k}(\tilde\beta(z))+\partial_{z_k}(\tilde\alpha(z)0=0,
\end{align*}
i.e.,
\begin{align}\label{sm.322.8}
\tilde\beta(z)\sideset{}{_{k=1}^m}{\sum}\partial_{z_k}(\tilde\alpha(z)0+\tilde\alpha(z)\sideset{}{_{k=1}^m}{\sum}\partial_{z_k}(\tilde\beta(z))+\sideset{}{_{k=1}^m}{\sum} \partial_{z_k}(\tilde\alpha(z))=0.
\end{align}
Using \eqref{sm.322.3} and \eqref{sm.322.4} to \eqref{sm.322.8}, we get
\begin{align}\label{sm.322.9} 
\tilde\beta^2(z)+\tilde\alpha(z)\sideset{}{_{k=1}^m}{\sum}\partial_{z_k}(\tilde\beta(z))+\tilde\beta(z)=0.
\end{align}
\medskip

\noindent
First suppose $\sum_{k=1}^m\partial_{z_k}(\tilde\beta(z)0\not\equiv 0$. Now in view of (\ref{sm.322.6}) and using Lemma \ref{L.3} to (\ref{sm.322.9}), we get  $\parallel T(r, \tilde\alpha)=o(T(r, \tilde\alpha))$ and so from Lemma \ref{L.7} and (\ref{sm.322.1}), we get 
$\parallel T(r,I_f)=o(T(r,I_f))$, which is impossible.
Next suppose $\sum_{k=1}^m\partial_{z_k}(\tilde\beta(z))\equiv 0$. Then from (\ref{sm.322.9}), we have
\begin{align}\label{sm.322.10} 
\tilde\beta^2(z)=-\tilde\beta(z).
\end{align}
If $\tilde\beta(z)$ is a transcendental entire function, then using Lemma \ref{L.4} to (\ref{sm.322.10}), we get 
$\parallel m(r,\tilde\beta)=o(T(r, \tilde\beta))$. This implies $\parallel T(r,\tilde\beta)=o(T(r,\tilde\beta))$, which is impossible. 
Hence $\tilde\beta(z)$ is a polynomial in $\mathbb{C}^m$. Therefore from (\ref{sm.322.10}), we find that $\tilde\beta(z)$ is a constant, say $k_1(\neq 0)$. Then from (\ref{sm.322.10}), we have 
$k_1^2=-k_1$, which implies that $k_1=-1$. Consequently
\begin{align}\label{sm.322.11} 
\tilde\beta(z)=\sideset{}{_{k=1}^m}{\sum}\partial_{z_k}(\alpha(z))=-1.
\end{align} 
Again from (\ref{sm.322.5}) and \eqref{sm.322.7}, we get 
\begin{align}\label{sm.322.12}
\sideset{}{_{k=1}^m}{\sum}\partial_{z_k}(P(z))=0.
\end{align}
The Lagrange's auxiliary equations corresponding to (\ref{sm.322.11}) are as follows
\begin{align}\label{sm.322.13} 
	\displaystyle
	\frac{d z_1}{1}=\frac{dz_2}{1}=\ldots=\frac{dz_m}{1}=\frac{d\tilde \alpha(z)}{-1}=\frac{dz_1+dz_2+\ldots+dz_m}{m}.
\end{align}
Clearly $z_2-z_1=d_1$, $z_3-z_2=d_2$, $\ldots$, $z_m-z_{m-1}=d_{m-1}$ and $m\tilde \alpha(z)+(z_1+z_2+\ldots+z_m)=d_m$ be independent solutions of (\ref{sm.322.11}). Therefore the general solution of (\ref{sm.322.11}) is as follows
\begin{align}\label{sm.322.14} 
	\displaystyle
	\tilde \alpha(z)=-(z_1+z_2+\ldots+z_m)/m+g_1(z),
\end{align}
where $g_1(z)=g_1(z_2-z_1,z_3-z_2,\ldots,z_m-z_{m-1})$ is an arbitrary finite-order transcendental entire function. Since $\sum_{j=1}^m\partial_{z_j}(f(z))$ is $c$-periodic, it follows from \eqref{sm.322.3}
\begin{align}\label{sm.322.15}	\displaystyle
\sideset{}{_{j=1}^m}{\sum}\partial_{z_j}(f(z)0=-(z_1+z_2+\ldots+z_m)/m+g_1(z).
\end{align}
Finally from \eqref{sm.322.1}, we have 
\begin{align*}	\displaystyle
\ol{f}(z)+f(z)=-\left(-(z_1+z_2+\ldots+z_m)/m+g_1(z)\right)^2+P(z),
\end{align*}
where $g_1(z)=g_1(z_2-z_1,z_3-z_2,\ldots,z_m-z_{m-1})$ is an arbitrary finite-order transcendental entire function such that (\ref{sm.322.15}) holds.
\end{proof}

\vspace{0.1in}
\noindent
{\bf Compliance of Ethical Standards:}\par

\noindent
{\bf Conflict of Interest.} The authors declare that there is no conflict of interest regarding the publication of this paper.\par

\noindent
{\bf Funding.} The authors have not received any funds for this research.\par

\noindent
{\bf Data availability statement.} Data sharing not applicable to this article as no data sets were generated or analysed during the current study.

\end{document}